\documentclass{article}

\usepackage{microtype}
\usepackage[nottoc]{tocbibind}
\usepackage[utf8]{inputenc}
\usepackage{graphicx}
\usepackage{float}
\usepackage{bm}
\usepackage{amsthm}
\usepackage{amsmath}
\usepackage{amssymb}
\usepackage{mathrsfs}
\usepackage{dsfont}
\usepackage[all]{xy}
\usepackage{tikz-cd}
\usepackage{stmaryrd}
\usepackage{enumitem}
\usepackage{tabularx}
\usepackage{placeins}
\usepackage{pdflscape}
\usepackage{hyperref}

\theoremstyle{definition}
\newtheorem{Definition}{Definition}[subsection]

\theoremstyle{plain}
\newtheorem{Theorem}[Definition]{Theorem}

\newcounter{mainthm}

\newtheorem{maintheorem}[mainthm]{Theorem}

\newtheorem{maincorollary}[mainthm]{Corollary}

\theoremstyle{plain}
\newtheorem{Proposition}[Definition]{Proposition}

\theoremstyle{plain}
\newtheorem{Lemma}[Definition]{Lemma}

\theoremstyle{plain}

\theoremstyle{plain}
\newtheorem{Corollary}[Definition]{Corollary}

\theoremstyle{plain}
\newtheorem{Conjecture}[Definition]{Conjecture}

\theoremstyle{plain}

\theoremstyle{definition}
\newtheorem{Construction}[Definition]{Construction}

\theoremstyle{definition}
\newtheorem{Example}[Definition]{Example}

\theoremstyle{definition}

\theoremstyle{remark}
\newtheorem{Remark}[Definition]{Remark}

\title{Compact Semisimple Tensor 2-Categories are Morita Connected}
\author{Thibault D.
Décoppet and Sean Sanford}
\date{May 2025}

\begin{document}

\bibliographystyle{alpha}

\maketitle
    \hspace{1cm}
    \begin{abstract}
    In \cite{D9}, it was shown that, over an algebraically closed field of characteristic zero, every fusion 2-category is Morita equivalent to a connected fusion 2-category, that is, one arising from a braided fusion 1-category.
    This result has recently allowed for a complete classification of fusion 2-categories.
    Here we establish that compact semisimple tensor 2-categories, which generalize fusion 2-categories to an arbitrary field of characteristic zero, also enjoy this ``Morita connectedness'' property.
    In order to do so, we generalize to an arbitrary field of characteristic zero many well-known results about braided fusion 1-categories over an algebraically closed field.
    Most notably, we prove that the Picard group of any braided fusion 1-category is indfinite, generalizing the classical fact that the Brauer group of a field is torsion.
    As an application of our main result, we derive the existence of braided fusion 1-categories indexed by the fourth Galois cohomology group of the absolute Galois group that represent interesting classes in the appropriate Witt groups.
    \end{abstract}

\section*{Introduction}

Throughout, we work over a fixed arbitrary field $\mathds{k}$ of characteristic zero.
Recall that a fusion 1-category over $\mathds{k}$ is a finite semisimple $\mathds{k}$-linear monoidal 1-category that is rigid, i.e.\ objects have duals, and has a simple monoidal unit.
This notion can be viewed as a wide generalization of that of a finite group.
More accurately, the definition of a fusion 1-category categorifies that of a separable algebra.
In more detail, it is well-known that separable algebras are precisely the fully dualizable objects in the sense of \cite{BD,L} of the symmetric monoidal Morita 2-category of $\mathds{k}$-algebras.
In fact, recall that a $\mathds{k}$-algebra is separable if and only if it is finite semisimple.
It was shown in \cite{DSPS13} that a fusion 1-category is a fully dualizable object of a certain symmetric monoidal Morita 3-category if and only if it separable, i.e.\ its Drinfeld center is finite semisimple.
But, it was established in \cite{ENO1} that, over an algebraically closed field, the Drinfeld center of a fusion 1-category is always finite semisimple.
It then follows that every fusion 1-category over $\mathds{k}$ is separable.
For completeness, let us also point out that a related perspective was given in \cite{JF} using the notion of higher condensation monad from \cite{GJF}.

Now, associated to the field $\mathds{k}$, there is a corresponding Brauer group $\mathrm{Br}(\mathds{k})$, whose objects are Morita equivalence classes of central simple algebras over $\mathds{k}$, or, equivalently, separable algebras whose center is identified with $\mathds{k}$.
In the context of higher categories, $\mathrm{Br}(\mathds{k})$ is the group of equivalence classes of invertible objects of the symmetric monoidal 2-category of $\mathds{k}$-algebras.
One of the key properties of the Brauer group is that, as shown by Noether \cite[Chapter II]{Jac}, it can be identified with the Galois cohomology group $H^2(\mathrm{Gal}(\overline{\mathds{k}}/\mathds{k}), \overline{\mathds{k}}^\times)$, where $\overline{\mathds{k}}$ denotes the algebraic closure of $\mathds{k}$.
Motivated by the analogy between separable algebras and fusion 1-categories recalled above, it is natural to consider the group of invertible objects of the symmetric monoidal Morita 3-category of fusion 1-categories over $\mathds{k}$.
We use $\mathrm{MoInv}(\mathbf{Vect}_\mathds{k})$ to denote this group.
Inspired by the classical identification of the Brauer group, it was shown in \cite{SS} that $\mathrm{MoInv}(\mathbf{Vect}_\mathds{k})$ can be identified with $H^3(\mathrm{Gal}(\overline{\mathds{k}}/\mathds{k}), \overline{\mathds{k}}^\times)$.

\begin{table}[!hbt]
\renewcommand{\arraystretch}{2.3}
    \begin{center}
    \caption{Some groups of interest}\label{tab:ourbelovedgroups}
    \begin{tabular}{|c|c|c|}
        \hline
        Name & Notation & Elements \\
        \hline
        Brauer group & $\mathrm{Br}(\mathds{k})$ & \parbox{4.5cm}{Morita equivalence classes of simple algebras with center $\mathds{k}$}\\
        \hline
         & $\mathrm{MoInv}(\mathbf{Vect}_{\mathds{k}})$ & \parbox{6cm}{Morita equivalence classes of fusion 1-categories with Drinfeld center $\mathbf{Vect}_{\mathds{k}}$} \\
         \hline
         Witt group & $\mathrm{Witt}(\mathbf{Vect}_{\mathds{k}})$ & \parbox{6.5cm}{Witt equivalence classes of braided fusion 1-categories with symmetric center $\mathbf{Vect}_{\mathds{k}}$}\\
         \hline
    \end{tabular}
    \end{center}
\end{table}

Working over $\overline{\mathds{k}}$, an algebraically closed field of characteristic zero, the notion of a fusion 2-category was introduced in \cite{DR}.
Categorifying the definition of a fusion 1-category over $\overline{\mathds{k}}$, a fusion 2-category is a rigid monoidal finite semisimple 2-category whose monoidal unit is simple.
This is based on the notion of finite semisimple 2-category introduced in \cite{DR}.
Succinctly, a finite semisimple 2-category is equivalent to the 2-category of finite semisimple module categories over a multifusion 1-category over $\overline{\mathds{k}}$.\footnote{Finite semisimple 2-categories admit a conceptual characterization \cite{DR}.}
An object of a finite semisimple 2-category is called simple if its multifusion 1-category of endomorphisms is fusion.
As observed in \cite{DR}, a key feature of finite semisimple 2-categories is that the Schur lemma fails:\ There can exist non-zero 1-morphisms between non-equivalent simple objects.
This leads to a notion of \textit{connected component}.
Two simple objects lie in the same connected component if there is a non-zero morphism bewteen them.
In particular, a coarse but useful invariant of a fusion 2-category is its set of connected components.
The most fundamental examples of fusion 2-categories are the connected fusion 2-categories, which are completely described by the braided fusion 1-category over $\overline{\mathds{k}}$ of endomorphisms of their monoidal unit.
However, fusion 2-categories are typically not connected. Many examples are discussed in \cite{DR} such as the fusion 2-categories of (twisted) 2-group-graded 2-vectors spaces, or the fusion 2-categories of 2-representations of a 2-group.

It was shown in \cite{D9} that every fusion 2-category over $\overline{\mathds{k}}$ is Morita equivalent to a connected one.
This was used in \cite{D10} in order to establish that fusion 2-categories over $\overline{\mathds{k}}$ are fully dualizable objects of an appropriate symmetric monoidal Morita 4-category.
It is then natural to consider the group of (Morita equivalence classes of) invertible objects of the symmetric monoidal Morita 4-category of fusion 2-categories over $\overline{\mathds{k}}$.
It was observed in \cite{D8} that this last group coincides with the Witt group $\mathrm{Witt}(\mathbf{Vect}_{\overline{\mathds{k}}})$ of non-degenerate braided fusion 1-categories over $\overline{\mathds{k}}$ introduced in \cite{DMNO}.
For completeness, we ought to mention that there is a symmetric monoidal Morita 4-category of braided fusion 1-categories over $\overline{\mathds{k}}$, and that $\mathrm{Witt}(\mathbf{Vect}_{\overline{\mathds{k}}})$ is isomorphic to its group of invertible objects \cite{BJSS}.
While the groups $\mathrm{Br}(\overline{\mathds{k}})$ and $\mathrm{MoInv}(\overline{\mathds{k}})$ are both trivial, the group $\mathrm{Witt}(\mathbf{Vect}_{\overline{\mathds{k}}})$ is infinite, and its precise structure has been the subject of much attention \cite{DMNO,DNO,NRWZ,JFR}.
More generally, for any symmetric fusion 1-category $\mathcal{E}$, there is a corresponding Witt group $\mathrm{Witt}(\mathcal{E})$ constructed from braided fusion 1-categories whose symmetric center is identified with $\mathcal{E}$, which was originally introduced in \cite{DNO}.
Using the fact that every fusion 2-category is Morita equivalent to a connected one, these Witt groups are intimately related to the Morita equivalence classes of fusion 2-categories.
This relation can be used to exhibit interesting elements in the groups $\mathrm{Witt}(\mathcal{E})$ as done in \cite{D9}, but can also be leveraged to classify fusion 2-categories over $\overline{\mathds{k}}$ up to monoidal equivalence \cite{DHJFNPPRY}.

Over the arbitrary field $\mathds{k}$, one can consider the notion of a compact semisimple tensor 2-category \cite{D5}.
Categorifying the definition of a fusion 1-category over $\mathds{k}$, a compact semisimple tensor 2-category is a rigid monoidal compact semisimple 2-category over $\mathds{k}$ whose monoidal unit is simple.
In particular, a compact semisimple 2-category is the 2-category of finite semisimple module categories over a multifusion 1-category over $\mathds{k}$.
The change in terminology from ``finite'' to ``compact'' is motivated by the following observation:\ While finite semisimple 2-categories over $\overline{\mathds{k}}$ have finitely many simple objects, compact semisimple 2-categories generally have infinitely many.
This can be seen very explicitly from the fact that the simple objects of the 2-category $\mathbf{2Vect}_{\mathds{k}}$ of finite semisimple $\mathds{k}$-linear 1-categories correspond to finite division $\mathds{k}$-algebras, which include in particular all finite field extensions of $\mathds{k}$.
On the other hand, compact semisimple 2-categories still have finitely many connected components, that is, equivalence classes of simple objects with non-zero 1-morphisms between them, and recover the notion of a finite semisimple 2-category over $\overline{\mathds{k}}$ from \cite{DR}.
The finite set of connected components of a compact semisimple tensor 2-category is again a useful invariant, and connected compact semisimple tensor 2-categories are still very closely related to braided fusion 1-category $\mathds{k}$.
We emphasize that a generic compact semisimple tensor 2-category has multiple connected components.

It was shown in \cite{D10} that every compact semisimple tensor 2-category over $\mathds{k}$ is a fully dualizable object of an appropriate symmetric monoidal Morita 4-category.
One may therefore think of compact semisimple tensor 2-categories over $\mathds{k}$ as a double categorification of separable $\mathds{k}$-algebras.
In particular, as a higher categorical version of the Brauer group, we can consider the group of (equivalence classes of) Morita invertible compact semisimple tensor 2-categories.
It is interesting to ask how this group relates to the Witt group $\mathrm{Witt}(\mathbf{Vect}_{\mathds{k}})$ of braided fusion 1-categories whose symmetric center is $\mathbf{Vect}_{\mathds{k}}$.
Furthermore, motivated by the lower categorical results, it is natural to wonder what is the precise relation between the last two groups and the Galois cohomology group $H^4(\mathrm{Gal}(\overline{\mathds{k}}/\mathds{k}), \overline{\mathds{k}}^\times)$.
More generally, it would be interesting to thoroughly investigate the structure of the groups $\mathrm{Witt}(\mathcal{E})$, for an arbitrary symmetric fusion 1-category $\mathcal{E}$ over $\mathds{k}$.

\subsection*{Results}

Our main theorem is the following generalization of \cite[Theorem 4.2.2]{D9} to arbitrary fields of characteristic zero.

\begin{maintheorem}\label{maintheorem}
Every compact semisimple tensor 2-category is Morita equivalent to a connected compact semisimple tensor 2-category.
\end{maintheorem}

\noindent As an immediate corollary, this implies that the group of Morita invertible compact semisimple tensor 2-categories and the Witt group $\mathrm{Witt}(\mathbf{Vect}_{\mathds{k}})$ of braided fusion 1-categories with symmetric center given by $\mathbf{Vect}_{\mathds{k}}$ are isomorphic. We will use this identification to exhibit non-trivial classes in $\mathrm{Witt}(\mathbf{Vect}_{\mathds{k}})$. By analogy with \cite{DHJFNPPRY}, we also expect that our main theorem can be used to classify compact semisimple tensor 2-categories. Though we relegate this second line of investigation to future work.

We now summarize the key steps of the proof.
Firstly, given any compact semisimple tensor 2-category $\mathfrak{C}$ over $\mathds{k}$, we can consider $\Omega\mathfrak{C}$, the braided fusion 1-category of endomorphisms of the monoidal unit.
The symmetric center $\mathcal{Z}_{(2)}(\Omega\mathfrak{C})$ is a symmetric fusion 1-category and therefore admits a fiber functor to the 1-category of super-vector spaces $\mathbf{SVect}_{\mathbb{K}}$ over some finite field extension $\mathbb{K}/\mathds{k}$ by Deligne's Theorem \cite{De}.
As in \cite{D9}, one can use this fiber functor to build a Morita equivalent compact semisimple tensor 2-category $\mathfrak{D}$ such that $\Omega\mathfrak{D}$ is either non-degenerate or slightly degenerate, meaning that $\mathcal{Z}_{(2)}(\Omega\mathfrak{C})$ is either $\mathbf{Vect}_{\mathbb{K}}$ or $\mathbf{SVect}_{\mathbb{K}}$.

Over an algebraically closed field $\overline{\mathds{k}}$ of characteristic zero, it was established in \cite{JFY} that if $\mathfrak{C}$ is a fusion 2-category over $\overline{\mathds{k}}$ with $\Omega\mathfrak{C}$ given by either $\mathbf{Vect}_{\overline{\mathds{k}}}$ or $\mathbf{SVect}_{\overline{\mathds{k}}}$, then every simple object of $\mathfrak{C}$ is invertible.
In \cite{D9}, this property was used to show that every such fusion 2-category $\mathfrak{C}$ over $\overline{\mathds{k}}$ is Morita equivalent to a connected fusion 2-category.
In fact, more generally, if $\mathfrak{C}$ is a fusion 2-category over $\overline{\mathds{k}}$ such that $\Omega\mathfrak{C}$ is either non-degenerate of slightly degenerate, then there exists an invertible object in every connected component of $\mathfrak{C}$ (see \cite[Corollary 4.3]{DHJFNPPRY} or Lemma \ref{lem:connectedcomponentsgroup} below).
These results do not hold over the arbitrary field $\mathds{k}$ (see for instance \cite[Remark 4.8]{D11}).
Nevertheless, for our present purposes, it is enough to establish that every compact semisimple tensor 2-category $\mathfrak{C}$ over $\mathds{k}$ with $\Omega\mathfrak{C}$ non-degenerate or slightly degenerate is Morita equivalent to a compact semisimple tensor 2-category that has an invertible object in every connected component.
At the decategorified level, that is, for fusion 1-category over $\mathds{k}$, a similar question was considered in \cite{SS} and solved using the procedure of inflation.
We develop the 2-categorical analogue, and use it to prove the next result.

\begin{maintheorem}
Let $\mathfrak{C}$ be a compact semisimple tensor 2-category with $\Omega\mathfrak{C}$ a non-degenerate, resp.\ slightly degenerate, braided fusion 1-category.
Then, $\mathfrak{C}$ is Morita equivalent to a compact semisimple tensor 2-category that has an invertible object in every connected component.
\end{maintheorem}

As in \cite{D9}, it will therefore be enough to show that any compact semisimple tensor 2-category $\mathfrak{C}$ that has an invertible object in every connected component is Morita equivalent to a connected one.
In order to do so, it is enough to construct a strongly connected rigid algebra $A$ in $\mathfrak{C}$ whose underlying object is a direct sum of invertible objects of $\mathfrak{C}$, with at least one in every connected component.
Note that the invertible objects appearing in such an algebra $A$ form a finite subgroup of the group of invertible objects of $\mathfrak{C}$.
While it is easy to construct finitely generated subgroups containing an invertible object in every connected component, there is a priori no guarantee that such subgroups are finite.
Namely, the group of invertible objects of a compact semisimple tensor 2-category is generally infinite.
This stems from the fact that, the group of invertible objects of a connected compact semisimple tensor 2-category, that is, the Picard group of the corresponding braided fusion 1-category, is infinite.
But, Picard groups generalize the classical Brauer groups of fields, so that it is natural to expect that they are torsion. In fact, as Brauer groups are also abelian, it follows that every finitely generated subgroup is finite.
We prove that Picard groups of braided fusion 1-categories over $\mathds{k}$ also satisfy this property.

\begin{maintheorem}
Let $\mathcal{B}$ be a braided fusion 1-category over $\mathds{k}$.
The Picard group $\mathrm{Pic}(\mathcal{B})$ is indfinite, i.e.\ every finitely generated subgroup is finite.
\end{maintheorem}

\begin{maincorollary}
Let $\mathfrak{C}$ be compact semisimple 2-category that has (at least) one invertible object in every connected component, then it is Morita equivalent to a connected compact semisimple tensor 2-category.
\end{maincorollary}

\noindent This concludes our survey of the proof of our main theorem.

\vspace{2mm}

As an application, we obtain the following result, which is essentially a categorification of \cite[Theorem 5.9]{SS}, and therefore a double categorification of the classical identification by Noether of the Brauer group with the second Galois cohomology of the absolute Galois group.

\begin{maintheorem}
There is a monomorphism $$H^4(\mathrm{Gal}(\overline{\mathds{k}}/\mathds{k});\overline{\mathds{k}}^\times)\hookrightarrow \mathrm{Witt}(\mathbf{Vect}_{\mathds{k}}).$$
\end{maintheorem}

\noindent More precisely, to every class in $H^4(\mathrm{Gal}(\overline{\mathds{k}}/\mathds{k});\overline{\mathds{k}}^\times)$, we associate a Morita equivalence class of invertible compact semisimple tensor 2-categories over $\mathds{k}$.
We show that this assignment is injective by using inflation.

Over the algebraically closed field $\overline{\mathds{k}}$ of characteristic zero, it was shown in \cite[Theorem 1.13]{DHJFNPPRY} that the assignment sending a symmetric fusion 1-category $\mathcal{E}$ to the corresponding Witt space is functorial.
Moreover, it is expected (see \cite[Theorem 1.14]{DHJFNPPRY}) that this functor preserves homotopy limits indexed by finite groups.
Thanks to our main theorem, we believe that these statements generalize to an arbitrary field $\mathds{k}$, and in Subsection \ref{sub:SpectralSeq} we discuss how our last theorem above can then be recovered from the corresponding homotopy fixed points spectral sequence.\footnote{The previous version of this work contained an appendix claiming to exhibit a
representative for the non-degenerate braided fusion 1-category over the real
numbers $\mathbb{R}$ corresponding to the non-trivial class in $H^4(\mathrm{Gal}(\mathbb{C}/\mathbb{R});\mathbb{C}^\times)\cong\mathbb{Z}/2$. The relevant computation turned out to be flawed. A correct representative is given in \cite[Theorem 5.9]{Sa2}.}

\subsubsection*{Acknowledgments}

We would like to thank Theo Johnson-Freyd for conversations related this topic, as well as all the organizers of the 2023 BIRS workshop on ``Subfactors and Fusion (2-)Categories (23w5091)" for providing a stimulating environment in the early stage of the project.
It is also our pleasure to thank Christopher Douglas for many helpful discussions.
The work of T.D.\ is supported in part by the Simons Collaboration on Global Categorical Symmetries.
The work of S.S. is partially supported by the
National Science Foundation under Grant No. DMS-2154389, for which he would specifically like to thank Dave Penneys.

\section{Preliminaries}

\subsection{Compact Semisimple Tensor 2-Categories}

Let $\mathds{k}$ be an arbitrary field of characteristic zero.
A multifusion 1-category over $\mathds{k}$ is a finite semisimple $\mathds{k}$-linear tensor 1-category, i.e.\ it carries a rigid monoidal structure.
A fusion 1-category is a multifusion 1-category whose monoidal unit $I$ is a simple object.
Given a multifusion 1-category $\mathcal{C}$ over $\mathds{k}$, we write $\Omega\mathcal{C} := End_{\mathcal{C}}(I)$, which is a direct sum of finite field extensions of $\mathds{k}$.
The fundamental properties of fusion 1-category over $\mathds{k}$ were investigated in \cite{Sa}.
Given any field extension $\mathbb{K}/\mathds{k}$, we define $\mathcal{C}_{\mathbb{K}}$, the base change of $\mathcal{C}$ along $\mathbb{K}/\mathds{k}$, to be the multifusion 1-category over $\mathbb{K}$ given by the Cauchy completion of $\mathcal{C}\otimes_{\mathds{k}}\mathbb{K}$.
In particular, we have $\Omega\mathcal{C}_{\mathbb{K}}=\Omega\mathcal{C}\otimes_{\mathds{k}}\mathbb{K}$, so that $\mathcal{C}_{\mathbb{K}}$ might be multifusion even if $\mathcal{C}$ is fusion.
When $\mathbb{K}=\overline{\mathds{k}}$ is the algebraic closure of $\mathds{k}$, we will also write $\overline{\mathcal{C}}:=\mathcal{C}_{\overline{\mathds{k}}}$.

Given a multifusion 1-category $\mathcal{C}$ over $\mathds{k}$, we will write $\mathbf{Mod}(\mathcal{C})$ for the corresponding 2-category of finite semisimple $\mathds{k}$-linear right $\mathcal{C}$-module 1-categories, $\mathds{k}$-linear $\mathcal{C}$-module functors, and $\mathcal{C}$-module natural transformations.
Such 2-categories are compact semisimple $\mathds{k}$-linear 2-categories as introduced in \cite{D5}.
More precisely, a $\mathds{k}$-linear 2-category $\mathfrak{C}$ is locally finite semisimple if $End_{\mathfrak{C}}(C)$ is a multifusion 1-category for every object $C$ of $\mathfrak{C}$.
An object $C$ of a locally finite semisimple 2-category $\mathfrak{C}$ is called simple if $End_{\mathfrak{C}}(C)$ is a fusion 1-category.
There is a notion of Cauchy completion for such 2-categories that was introduced in \cite{GJF,DR} (see also \cite{D1}).
If we assume in addition that $\mathfrak{C}$ is Cauchy complete, then every object splits as a finite direct sum of simple objects.
We say that two simple objects are connected if there exists a non-zero 1-morphism between them.
This defines an equivalence relation on the set of simple objects, and we write $\pi_0(\mathfrak{C})$ for the corresponding set of equivalence classes, which we refer to as the connected components of $\mathfrak{C}$.
\begin{Definition}
A compact semisimple $\mathds{k}$-linear 2-category is a Cauchy complete locally finite $\mathds{k}$-linear 2-category that has finitely many connected components.
\end{Definition}
\noindent It was shown in \cite{D5} that every compact semisimple tensor 2-category is of the form $\mathbf{Mod}(\mathcal{C})$ for some multifusion 1-category $\mathcal{C}$ over $\mathds{k}$.

When $\mathds{k}=\overline{\mathds{k}}$ is algebraically closed, every connected component contains finitely many (equivalence classes of) simple objects.
In that case, the notion of a compact semisimple 2-category recovers that of a finite semisimple 2-category introduced in \cite{DR}.
Slightly more generally, it was shown in \cite{D5} that a compact semisimple 2-category has finitely many simple objects if and only if $\mathds{k}$ is algebraically closed or real closed.

\begin{Definition}
A compact semisimple multitensor 2-category over $\mathds{k}$ is a compact semisimple 2-category equipped with a rigid monoidal structure, that is every object has a right and left dual.
A compact semisimple tensor 2-category is a compact semisimple multitensor 2-category with simple monoidal unit $I$.
\end{Definition}
\noindent For any compact semisimple multitensor 2-category $\mathfrak{C}$, we write $\Omega\mathfrak{C}:=End_{\mathfrak{C}}(I)$ for the braided multifusion 1-category of endomorphisms of the monoidal unit.
By definition, $\Omega\mathfrak{C}$ is a fusion 1-category if and only if $\mathfrak{C}$ is a compact semisimple tensor 2-category.
Every connected compact semisimple tensor 2-category is of the form $\mathbf{Mod}(\mathcal{B})$ for a braided fusion 1-category $\mathcal{B}$ over $\mathds{k}$.
Moreover, if $\mathfrak{C}$ is a compact semisimple tensor 2-category, then $\mathfrak{C}^0$, the connected component of the identity of $\mathfrak{C}$, satisfies $\mathfrak{C}^0\simeq\mathbf{Mod}(\mathcal{B})$.
Over an algebraically closed field $\mathds{k}=\overline{\mathds{k}}$, this notion of compact semisimple tensor 2-category is exactly that of a fusion 2-category from \cite{DR}.

Given any field extension $\mathbb{K}/\mathds{k}$ and compact semisimple 2-category $\mathfrak{C}$ over $\mathds{k}$, we define $\mathfrak{C}_{\mathbb{K}}$, the base change of $\mathfrak{C}$ along $\mathbb{K}/\mathds{k}$, to be the compact semisimple 2-category over $\mathbb{K}$ given by the Cauchy completion of the local Cauchy completed tensor product $\mathfrak{C}\widehat{\otimes}_{\mathds{k}}\mathbb{K}$.
In particular, it follows from the definitions that if $\mathfrak{C}\simeq \mathbf{Mod}(\mathcal{C})$ as compact semisimple 2-categories over $\mathds{k}$ for some multifusion 1-category $\mathcal{C}$ over $\mathds{k}$, then $\mathfrak{C}_{\mathbb{K}}\simeq \mathbf{Mod}(\mathcal{C}_{\mathbb{K}})$ as compact semisimple 2-categories over $\mathbb{K}$.
When $\mathbb{K}=\overline{\mathds{k}}$ is the algebraic closure of $\mathds{k}$, we will also write $\overline{\mathfrak{C}}:=\mathfrak{C}_{\overline{\mathds{k}}}$.
If $\mathfrak{C}$ is a compact semisimple multitensor 2-category over $\mathds{k}$, then its base change $\mathfrak{C}_{\mathbb{K}}$ is a compact semisimple multitensor 2-category over $\mathbb{K}$.
For later use, let us record the fact that there is a canonical monoidal $\mathds{k}$-linear 2-functor $\mathfrak{C}\rightarrow \mathfrak{C}_{\mathbb{K}}$.
Finally, let us also observe that it follows from the definitions that $\Omega(\mathfrak{C}_{\mathbb{K}})\simeq (\Omega\mathfrak{C})_{\mathbb{K}}$, so that we may simply write $\Omega\mathfrak{C}_{\mathbb{K}}$.

\begin{Example}\label{ex:gradedCSS2C}
Let $\mathbb{K}/\mathds{k}$ be a finite Galois extension.
Let $G$ be a finite group acting on $\mathbb{K}$, i.e.\ there is a group homomorphism $f:G\rightarrow \mathrm{Gal}(\mathbb{K}/\mathds{k})$.
Associated to this data, there is a compact semisimple tensor 2-category $\mathbf{2Vect}_{\mathbb{K}}(G)$ that is linear over $\mathds{k}$.
The most straightforward way to construct this compact semisimple tensor 2-category is as the semidirect product $\mathbf{2Vect}_{\mathbb{K}}\rtimes_f G$.
Equivalently, it can be obtained by Cauchy completing the monoidal 2-category $\mathrm{B}^2\mathbb{K}\rtimes_f G$ given by the semi-direct product of the monoidal 2-category $\mathrm{B}^2\mathbb{K}$ that has only one object, one 1-morphism, and abelian group of 2-morphisms given by $\mathbb{K}$.
Slightly more generally, given a 4-cocycle representing a class $\pi$ in $H^4(G;\mathbb{K}^{\times})$, where $\mathbb{K}^{\times}$ is viewed as a $G$-module in the obvious way, we can define $\mathbf{2Vect}^{\pi}_{\mathbb{K}}(G)$ by using the 4-cocycle representing $\pi$ to twist the pentagonator 2-natural isomorphism.
This construction only depends on the cohomology class $\pi$.
Let us observe that the compact semisimple tensor 2-category $\mathbf{2Vect}^{\pi}_{\mathbb{K}}(G)$ is in fact linear over $\mathbb{K}^G$, the subfield of $\mathbb{K}$ fixed under the action of $G$.
\end{Example}

\subsection{Rigid Algebras, (Bi)Modules, and Morita Equivalences}

Let $\mathfrak{C}$ be a compact multitensor 2-category over $\mathds{k}$.
An algebra $A$ in $\mathfrak{C}$ is rigid if the multiplication map $A\Box A\rightarrow A$ admits a right adjoint $m^*$ as a 1-morphism of $A$-$A$-bimodules.
It was proven in \cite[Theorem 3.1.4]{D10} that, under our assumption that $\mathds{k}$ has characteristic zero, every such rigid algebra is separable, in the sense that the $A$-$A$-bimodules 2-morphism $m\circ m^*\Rightarrow Id_A$ admits a section.
We will use this fact repeatedly below.
In particular, if $A$ is a rigid algebra in $\mathfrak{C}$, the 2-categories $\mathbf{Mod}_{\mathfrak{C}}(A)$ of right $A$-module in $\mathfrak{C}$ and $\mathbf{Bimod}_{\mathfrak{C}}(A)$ of $A$-$A$-bimodules in $\mathfrak{C}$ are compact semisimple.
Moreover, $\mathbf{Bimod}_{\mathfrak{C}}(A)$ is a compact semisimple multitensor 2-cateogry \cite{D8}.
If $A$ is in addition connected, in the sense that the unit 1-morphism $u:I\rightarrow A$ is simple, then $\mathbf{Bimod}_{\mathfrak{C}}(A)$ is a compact semisimple tensor 2-category.

Let $\mathcal{B}$ be a braided fusion 1-category over $\mathds{k}$.
Connected rigid algebras in the compact semisimple tensor 2-category $\mathbf{Mod}(\mathcal{B})$ are precisely $\mathcal{B}$-central fusion 1-categories over $\mathds{k}$.
Recall that a $\mathcal{B}$-central fusion 1-category is a fusion 1-category $\mathcal{C}$ equipped with a braided tensor functor $\mathcal{B}\rightarrow \mathcal{Z}(\mathcal{C})$ to the Drinfeld center of $\mathcal{C}$.
The corresponding compact semisimple 2-category of modules is $\mathbf{Mod}(\mathcal{C})$, which has a natural left action by $\mathbf{Mod}(\mathcal{B})$.

We will now briefly discuss the concept of Morita equivalence between compact semisimple tensor 2-categories, referring the reader to \cite{D8} for more details.
In order to formulate the definition, recall that an algebra $A$ in $\mathfrak{C}$ is faithful if, for every simple summand $J$ of the monoidal unit $I$ of $\mathfrak{C}$, $A\Box J$ is nonzero.
In particular, if $\mathfrak{C}$ has simple monoidal unit, i.e.\ $\mathfrak{C}$ is a compact semisimple tensor 2-category, then every algebra is automatically faithful.
Two compact semisimple multitensor 2-categories $\mathfrak{C}$ and $\mathfrak{D}$ over $\mathds{k}$ are Morita equivalent if there exists a faithful rigid algebra $A$ in $\mathfrak{C}$ and an equivalence $\mathbf{Bimod}_{\mathfrak{C}}(A)\simeq \mathfrak{D}$ of monoidal 2-categories over $\mathds{k}$.
This defines an equivalence relation on compact semisimple multitensor 2-categories.
Moreover, this definition of Morita equivalence between compact semisimple multitensor 2-categories may also be reformulated using the language of module 2-categories, for which we refer the reader to \cite{D4} for the precise definitions.
Given a compact semisimple left $\mathfrak{C}$-module 2-category $\mathfrak{M}$, we write $\mathbf{End}_{\mathfrak{C}}(\mathfrak{M})$ for the compact semisimple tensor 2-category of left $\mathfrak{C}$-module endofunctors on $\mathfrak{M}$.
Then, $\mathfrak{C}$ and $\mathfrak{D}$ are Morita equivalent if (and only if) there exists a faithful compact semisimple left $\mathfrak{C}$-module 2-category $\mathfrak{M}$, and a monoidal equivalence $\mathfrak{D}\simeq \mathbf{End}_{\mathfrak{C}}(\mathfrak{M})^{\mathrm{mop}}$.
As is explained in \cite{D8}, that these two definitions coincide is a consequence of the fact that every compact semisimple left $\mathfrak{C}$-module is of the form $\mathbf{Mod}_{\mathfrak{C}}(A)$ for some rigid algebra $A$ in $\mathfrak{C}$ by \cite{D4}.
With $\mathds{k} = \overline{\mathds{k}}$ algebraically closed, Morita equivalence classes of fusion 2-categories were extensively studied in \cite{D9}.
In particular, it was shown therein that every fusion 2-category is Morita equivalent to a connected one.
We will prove that this statement holds more generally for any compact semisimple tensor 2-category over an arbitrary field (of characteristic zero).
In order to do so, we will make use of the following technical lemma, whose proof follows that of \cite[Lemma 4.2.1]{D9} verbatim. We say that an algebra $A$ is strongly connected if its unit 1-morphism $u:I\rightarrow A$ is the inclusion of a summand (as an object of $\mathfrak{C}$).

\begin{Lemma}\label{lem:EnoughInvertibleMoritaConnected}
Let $\mathds{k}$ be an arbitrary field of characteristic zero, and let $\mathfrak{C}$ be a compact semisimple tensor 2-category.
If $A$ is a strongly connected rigid algebra in $\mathfrak{C}$, whose underlying object contains an invertible object in every connected component of $\mathfrak{C}$, then $\mathbf{Bimod}_{\mathfrak{C}}(A)$ is a connected compact semisimple tensor 2-category.
\end{Lemma}

Given a compact semisimple tensor 2-category $\mathfrak{C}$ over $\mathds{k}$, or, more generally any monoidal 2-category, we use $\mathscr{Z}(\mathfrak{C})$ to denote its Drinfeld center.
The Drinfeld center is a canonical braided monoidal 2-category associated to $\mathfrak{C}$, and was introduced in \cite{Cr}.
It follows from \cite{D10} that the Drinfeld center of a compact semisimple (multi)tensor 2-category $\mathfrak{C}$ is a compact semisimple (multi)tensor 2-category.
For later use, let us record that $\mathscr{Z}(\mathfrak{C})$ is equivalent to the monoidal 2-category of $\mathfrak{C}$-$\mathfrak{C}$-bimodule endo-2-functors on $\mathfrak{C}$.
The next theorem, which categorifies a well-know result in the theory of fusion 1-categories, follows from the results of \cite[Section 2]{D9}.
Namely, even though these are stated over an algebraically closed field of characteristic zero, their proofs do not make use of this assumption.

\begin{Theorem}[{\cite[Thm.
2.3.2]{D9}}]\label{thm:MoritaCenter}
Let $\mathfrak{C}$ and $\mathfrak{D}$ be two Morita equivalent compact semisimple tensor 2-categories over $\mathds{k}$.
Then, there is an equivalence $\mathscr{Z}(\mathfrak{C})\simeq \mathscr{Z}(\mathfrak{D})$ of braided compact semisimple tensor 2-categories.
\end{Theorem}

\subsection{Galois Non-triviality and Galois Grading}\label{sub:Galois}

We now present a categorification of the theory of Galois non-triviality and the associated Galois grading developed in \cite{Sa,SS} for fusion 1-categories.

Given a compact semisimple tensor 2-category $\mathfrak{C}$ over $\mathds{k}$, the field $\mathbb{K}:=\Omega^2\mathfrak{C}$ plays a key role.
More precisely, $\mathbb{K}$ is a finite field extension of $\mathds{k}$, which is generically non-trivial, and may therefore interact in interesting ways with the objects of $\mathfrak{C}$.
More precisely, associated to any object $C$ of $\mathfrak{C}$ there are two braided functors $L_C:\Omega\mathfrak{C}\hookrightarrow\mathcal{Z}(End_{\mathfrak{C}}(C))$ and $R_C:(\Omega\mathfrak{C})^{\mathrm{rev}}\hookrightarrow\mathcal{Z}(End_{\mathfrak{C}}(C))$, that correspond to the left and right actions of $\Omega\mathfrak{C}$ on $End_{\mathfrak{C}}(C)$.
Then, upon taking loops, these restrict to two embeddings $\lambda_C,\rho_C:\Omega^2\mathfrak{C}\hookrightarrow\Omega\mathcal{Z}(End_{\mathfrak{C}}(C))$, called the left and right embeddings for $C$ respectively.

\begin{Definition}
The object $C$ of $\mathfrak{C}$ is said to be Galois trivial if $\lambda_C=\rho_C$.
Otherwise $C$ is said to be Galois non-trivial.
\end{Definition}

\noindent Note that $\Omega\mathcal{Z}(End_{\mathfrak C}(C))$ is always a commutative algebra, and is thus a product of finite field extensions of $\mathbb{K}$.
In particular, when $C$ is simple, $\lambda_C$ and $\rho_C$ are embeddings of fields.
Moreover, the field $\Omega\mathcal{Z}(End_{\mathfrak{C}}(C))$ is determined up to canonical isomorphism by the component $[C]\in\pi_0(\mathfrak{C})$ of $C$.
When isomorphisms are not a cause of concern, we will occasionally abuse notation by saying that $\Omega\mathcal{Z}(End_{\mathfrak{C}}(C))$ is a field.
The next lemma establishes that the isomorphism classes of $\lambda_C$ and $\rho_C$ do not depend on the choice of simple object $C$.

\begin{Lemma}
Let $C$ and $D$ be two simple objects in a fixed connected component of $\mathfrak{C}$.
There is a canonical isomorphism $\tau_{D,C}:\Omega\mathcal{Z}(End_{\mathfrak{C}}(C))\to\Omega\mathcal{Z}(End_{\mathfrak{C}}(D))$ that satisfies the following two equations
    \begin{gather}
        \tau_{D,C}\circ\lambda_C\;=\;\lambda_D\,,\hspace{4mm}\tau_{D,C}\circ\rho_C\;=\;\rho_D\,.\label{Eqn:TauLambdaRho}
    \end{gather}
\end{Lemma}

\begin{proof}
The finite semisimple bimodule 1-category $Hom_{\mathfrak{C}}(C,D)$ is a Morita equivalence between $End_{\mathfrak{C}}(C)$ and $End_{\mathfrak{C}}(D)$, and therefore gives rise to a braided equivalence $\mathcal{Z}(End_{\mathfrak{C}}(C))\to\mathcal{Z}(End_{\mathfrak{C}}(D))$, which restricts to an isomorphism of fields $\tau_{D,C}:\Omega\mathcal{Z}(End_{\mathfrak{C}}(C))\to\Omega\mathcal{Z}(End_{\mathfrak{C}}(D))$.
The desired equations follow.
\end{proof}

\begin{Example}\label{Eg:ComponentsOfBim(Vec_K)}
For any finite field extension $\mathbb{K}/\mathds{k}$, $\mathbf{Vect}_{\mathbb{K}}$ is a connected rigid algebra in $\mathbf{2Vect}_{\mathds{k}}$, so we can consider the compact semisimple tensor 2-category $\mathfrak{C}:=\mathbf{Bimod}_{\mathbf{2Vect}_{\mathds{k}}}(\mathbf{Vect}_{\mathbb{K}})$.
It follows by inspection that $\Omega\mathfrak{C}\simeq \mathbf{Vect}_{\mathbb{K}}$.
We also have that the underlying compact semisimple 2-category of $\mathfrak{C}$ is given by $\mathbf{Mod}(\mathbf{Vect}_{\mathbb{K}}\boxtimes_{\mathds{k}}\mathbf{Vect}_{\mathbb{K}})$.
We can identify $\mathbf{Vect}_{\mathbb{K}}\boxtimes_{\mathds{k}}\mathbf{Vect}_{\mathbb{K}}$ with the multifusion 1-category of finite dimensional modules over the finite semisimple commutative algebra $\mathbb{K}\otimes_{\mathds{k}}\mathbb{K}$.
Now, by the primitive element theorem, there exists $\theta\in\mathbb{K}$ so that $\mathbb{K}=\mathds{k}(\theta)$.
Let $f(x)\in\mathds{k}[x]$ be the minimal polynomial for $\theta$, and suppose $f(x)=\prod_{i=1}^nf_{i}(x)$ is a factorization of $f(x)$ over $\mathbb K[x]$ into irreducible polynomials.
Writing $\mathbb{K}_i:=\mathbb{K}[x]/(f_i(x))$, we therefore have $\mathbb{K}\otimes_{\mathds{k}}\mathbb{K}\cong \mathbb{K}[x]/(f(x))\cong \prod_{i=1}^n\mathbb{K}_i$.
Then, the connected components of $\mathfrak{C}$ are $\mathbf{2Vect}_{\mathbb{K}_i}$ generated by $\mathbf{Vect}_{\mathbb{K}_i}$, and $\Omega\mathcal{Z}(End_{\mathfrak{C}}(\mathbf{Vect}_{\mathbb{K}_i}))=\mathbb{K}_i$ for $i=1,...,n$.
Without loss of generality, we may assume that $f_1(x)=x-\theta$, whence $\mathbb{K}_1\cong\mathbb K$ corresponds to the connected component of the monoidal unit.
In the event that the extension $\mathbb{K}/\mathds{k}$ is Galois, every field $\mathbb{K}_i$ is isomorphic to $\mathbb{K}$.
\end{Example}

When $\mathbb K/\mathds k$ is not normal, the example above shows that the left and right embeddings are in general not isomorphisms.
Just as in the decategorified setting, it is useful to single out simple objects for which the left and right embeddings are isomorphisms.

\begin{Definition}\label{def:normal}
A simple object $C$ of $\mathfrak{C}$ is said to be normal if $\Omega\mathcal{Z}(End_{\mathfrak{C}}(C))\cong\Omega^2\mathfrak{C}=\mathbb{K}$, and a component $[C]$ of $\mathfrak{C}$ is said to be normal if $\Omega\mathcal{Z}(End_{\mathfrak{C}}(C))\cong\Omega^2\mathfrak{C}=\mathbb{K}$.
The compact semisimple tensor 2-category is normal if all of its connected components are normal.
\end{Definition}

\noindent In particular, equation \eqref{Eqn:TauLambdaRho} has the following consequence.

\begin{Lemma}\label{lem:DefOfg}
Suppose that $[C]$ is a normal connected component of $\mathfrak{C}$.
To any simple object $D\in [C]$, there is a well-defined automorphism $g_D:=\lambda^{-1}_D\circ\rho_D$ of $\mathbb{K}/\mathds{k}$.
This automorphism is independent of the choice of $D$ and so we define $g_{[C]}:=g_D$.
\end{Lemma}

\noindent Furthermore, categorifying Lemma 3.11 of \cite{PSS}, we have the following technical result.

\begin{Lemma}\label{lem:Galoisgradingtechnical}
Suppose $C$, $D$, and $E$ are normal simple objects of $\mathfrak{C}$.
If there exists a nonzero 1-morphism $f:C\Box D\to E$, then $g_E=g_C\circ g_D$.
\end{Lemma}
\begin{proof}
Let $\lambda\in\mathbb K=\Omega^2\mathfrak C$ be a scalar, and assume without loss of generality that $f$ is simple.
Using functoriality of $\Box$, the 2-morphism $f\Box \lambda$ is equal to both $(\mathrm{id}_E\Box \lambda)\circ f$ and $f\circ(\mathrm{id}_{C\Box D}\Box \lambda)$.
Using the definition of the automorphism $g$ and functoriality again, these are equal to $g_E(\lambda)\Box f$ and $(g_C\circ g_D)(\lambda)\Box f$ respectively.
Thus, we find that $$\big(g_E(\lambda)-(g_C\circ g_D)(\lambda)\big)\Box f=0$$ for any $\lambda$, but then simplicity of $f$ implies that $g_E(\lambda)-(g_C\circ g_D)(\lambda)=0$, and the claim follows.
\end{proof}

\begin{Proposition}\label{prop:Galoisgrading}
If all the connected components of $\mathfrak{C}$ are normal, then $\mathfrak{C}$ admits a (not necessarily faithful) grading by the group $\mathrm{Gal}(\mathbb K/\mathds k)$, where $\mathfrak{C}_g$ is spanned by all those simple objects $C$ for which $g_C=g$.
Moreover, the connected component $[C]$ of $\mathfrak{C}$ lies entirely within the graded component $\mathfrak{C}_{g_{[C]}}$.
\end{Proposition}
\begin{proof}
The first part follows from Lemma \ref{lem:Galoisgradingtechnical} and the second part follows from Lemma \ref{lem:DefOfg}.
\end{proof}

\begin{Example}\label{ex:MoritatrivialCSS2Cs}
Let $\mathbb{K}/\mathds{k}$ be a finite Galois extension, and write $\Gamma := \mathrm{Gal}(\mathbb{K}/\mathds{k})$.
For any class $\pi$ in $H^4(\Gamma;\mathbb{K}^\times)$, the compact semisimple tensor 2-category $\mathbf{2Vect}^{\pi}_{\mathbb{K}}(\Gamma)$ introduced in Example \ref{ex:gradedCSS2C} admits a faithful Galois grading.
Moreover, it acts on $\mathbf{2Vect}_{\mathbb{K}}$ if and only if $\pi=\mathrm{triv}$.
By direct computation, we find that the compact semisimple tensor 2-category of left $\mathbf{2Vect}_{\mathds{k}}$-module endo-2-functors of $\mathbf{2Vect}_{\mathbb{K}}$ is given by $\mathbf{2Vect}_{\mathbb{K}}(\Gamma)$.
Expressed differently, this means that there is an equivalence $$\mathbf{Bimod}_{\mathbf{2Vect}_{\mathds{k}}}(\mathbf{Vect}_{\mathbb{K}})\simeq\mathbf{2Vect}_{\mathbb{K}}(\Gamma)$$ of compact semisimple tensor 2-categories, thereby explicitly identifying the compact semisimple tensor 2-categories introduced in Example \ref{Eg:ComponentsOfBim(Vec_K)} when $\mathbb{K}/\mathds{k}$ is a finite Galois extension.
\end{Example}

\subsection{Group Graded Extensions}

Group graded compact semisimple tensor 2-categories will feature prominently.
This is motivated by the fact that group gradings play key roles both in the decategorified setting \cite{SS} and in the algebraically closed case \cite{D9}.
Group graded extensions of fusion 2-categories were studied extensively in \cite{D11}, where a categorification of the main result of \cite{ENO2} was obtained.
For our present purposes, we only need to recall the following result.

\begin{Proposition}[{\cite[Prop.\ 2.4 \& Rem.\ 2.13]{D11}}]
Let the compact semisimple tensor 2-category $\mathfrak{C} = \boxplus_{g\in G}\mathfrak{C}_g$ be a faithfully $G$-graded extension of $\mathfrak{C}_e$.
Then, for any $g\in G$, the $\mathfrak{C}_e$-$\mathfrak{C}_e$-bimodule 2-category $\mathfrak{C}_g$ is invertible, i.e.\ it witnesses a Morita autoequivalence of $\mathfrak{C}_e$.
\end{Proposition}

In general, group graded extensions of compact semisimple tensor 2-categories can be complicated (see \cite{DY23} and \cite[Remark 3.8]{D11}).
The more restrictive notion of a quasi-trivial extension is substantially more well-behaved.
This is well-known in the context of extensions of fusion 1-categories over an algebraically closed field of characteristic zero \cite[Section 4.3]{ENO2}.

\begin{Definition}
Let $G$ be a finite group, and let $\mathfrak{C}=\boxplus_{g\in G}\mathfrak{C}_g$ be a $G$-graded compact semisimple tensor 2-category.
We say that $\mathfrak{C}$ is a quasi-trivial extension if every graded component $\mathfrak{C}_g$ of $\mathfrak{C}$ is equivalent to $\mathfrak{C}_e$ as a left, or equivalently right, $\mathfrak{C}_e$-module 2-category.
\end{Definition}

\noindent The following lemma categorifies \cite[Lemma 3.1]{Gal}.

\begin{Lemma}\label{lem:quasitrivialinvertible}
Let $\mathfrak{C}$ be a faithfully $G$-graded compact semisimple tensor 2-category.
The extension $\mathfrak{C}=\boxplus_{g\in G}\mathfrak{C}_g$ is quasi-trivial if and only if every graded component $\mathfrak{C}_g$ of $\mathfrak{C}$ contains an invertible object.
\end{Lemma}
\begin{proof}
The backwards direction is clear, so let $\mathfrak{C}=\boxplus_{g\in G}\mathfrak{C}_g$ be quasi-trivial.
For any $g\in G$, there is an equivalence $\mathfrak{C}_g\simeq \mathfrak{C}_e$ of left $\mathfrak{C}_e$-module 2-categories.
But, by \cite[Corollary 2.5 and Remark 2.13]{D11}, there is an equivalence of $\mathfrak{C}_e$-$\mathfrak{C}_e$-bimodule 2-categories $\mathfrak{C}_{g^{-1}}\simeq \mathbf{Fun}_{\mathfrak{C}_e}(\mathfrak{C}_g,\mathfrak{C}_e)$, and we therefore find that $\mathfrak{C}_{g^{-1}}\simeq \mathfrak{C}_e$ as right $\mathfrak{C}_e$-module 2-categories.
Under these identifications, the 2-functor $\Box:\mathfrak{C}_{g^{-1}}\times \mathfrak{C}_g\rightarrow \mathfrak{C}_e$ is identified with $\Box:\mathfrak{C}_{e}\times \mathfrak{C}_e\rightarrow \mathfrak{C}_e$, so that the claim follows.
\end{proof}

The following technical results will be used subsequently.
Firstly, working over an algebraically closed field $\mathds{k}=\overline{\mathds{k}}$, recall that a strongly fusion 2-category is a fusion 2-category $\mathfrak{C}$ such that $\Omega\mathfrak{C}$ is equivalent to either $\mathbf{Vect}$ or $\mathbf{SVect}$ as a braided fusion 1-category.
It was shown in \cite{JFY} that every simple object of a strongly fusion 2-category is invertible. As a consequence, the set of connected components of such a fusion 2-category inherits a group structure from the monoidal product. This last fact generalizes to any fusion 2-category $\mathfrak{C}$ whose loops $\Omega\mathfrak{C}$ are either a non-degenerate or a slightly degenerate braided fusion 1-category, that is, the symmetric center of the braided fusion 1-category $\Omega\mathfrak{C}$ is either $\mathbf{Vect}$ or $\mathbf{SVect}$.

\begin{Lemma}[{\cite[Corollary 4.3]{DHJFNPPRY}}]\label{lem:connectedcomponentsgroup}
Let $\mathds{k}=\overline{\mathds{k}}$ be algebraically closed, and let $\mathfrak{C}$ be a fusion 2-category such that $\mathcal{B}:=\Omega\mathfrak{C}$ is a non-degenerate or slightly degenerate braided fusion 1-category.
Then, the set of connected components of $\mathfrak{C}$ inherits a group structure from the monoidal product of $\mathfrak{C}$, so that $\mathfrak{C}$ is a group graded extension of $\mathfrak{C}^0$.
\end{Lemma}

\noindent Secondly, we will also make use of the following related observation. Recall that a braided fusion 1-category $\mathcal{B}$ over $\overline{\mathds{k}}$ is called completely anisotropic if every commutative separable $A$ in $\mathcal{B}$ is a direct sum of copies of the monoidal unit. In particular, every completely anisotropic braided fusion 1-category must be either non-degenerate or slightly degenerate.

\begin{Lemma}\label{lem:Moritageneralizedstronglyfusion}
Let $\mathds{k}=\overline{\mathds{k}}$ be algebraically closed, and let $\mathfrak{C}$ be a fusion 2-category such that $\mathcal{B}:=\Omega\mathfrak{C}$ is a completely anisotropic braided fusion 1-category.
If $A$ is a connected rigid algebra in $\mathfrak{C}$ such that $\Omega\mathbf{Bimod}_{\mathfrak{C}}(A)$ is completely anisotropic, then $A$ lies in $\mathfrak{C}^0$.
Moreover, in that case, the forgetful 2-functors $\mathbf{Bimod}_{\mathfrak{C}}(A)\rightarrow\mathfrak{C}$ and $\mathbf{Mod}_{\mathfrak{C}}(A)\rightarrow\mathfrak{C}$ induce bijections on connected components.
\end{Lemma}
\begin{proof}
We claim that every connected algebra in such a fusion 2-category $\mathfrak{C}$ is necessarily strongly connected.
Namely, given a connected algebra $A$ in $\mathfrak{C}$, we find that the algebra $A^0$ given by the restriction of $A$ to the connected component of the identity $\mathfrak{C}^0$ is a $\mathcal{B}$-central finite semisimple monoidal 1-category $\mathcal{A}^0$.
By \cite[Lemma 5.12]{DMNO} and its obvious generalization to the slightly degenerate case, it follows that the monoidal functor $\mathcal{B}\rightarrow\mathcal{A}^0$ is fully faithful.
This is exactly the statement that the algebra $A^0$, and therefore also $A$ is strongly connected.

Let us now assume that $A$ is a connected rigid algebra in $\mathfrak{C}$, which is therefore automatically strongly connected.
The connected components of $\mathfrak{C}$ inherit a group structure, and the support of the rigid algebra $A$ is necessarily a subgroup by \cite[Lemma 4.10]{DHJFNPPRY}.
Up to replacing $\mathfrak{C}$ by a full fusion sub-2-category, we may therefore assume that $A$ contains a simple object in every connected component of $\mathfrak{C}$.
The first part of the statement is then a consequence of the fact that every fusion 2-category that is Morita equivalent to $\mathbf{Mod}(\mathcal{B})$ with $\mathcal{B}$ completely anisotropic must itself be connected by \cite[Corollary 3.1.5]{D9} and \cite[Corollary 3.26]{DMNO}.
The second part then follows from \cite[Lemma 4.10]{DHJFNPPRY}.
\end{proof}

\section{Braided Fusion 1-Categories over an Arbitrary Field}

We establish various properties of braided fusion 1-categories over an arbitrary field $\mathds{k}$ of characteristic zero.

\subsection{Centers and Non-Degeneracy}\label{sub:CenterNonDegeneracy}

Let $\mathcal{B}$ be a braided fusion 1-category over $\mathds{k}$ with braiding $c_{B,C}:B\otimes C\cong C\otimes B$.
We use $\mathcal{B}^{\mathrm{rev}}$ to denote fusion 1-category $\mathcal{B}$ equipped with the braiding $c_{C,B}^{-1}:B\otimes C\cong C\otimes B$.
We begin by recording the following fact, which will be used repeatedly.

\begin{Lemma}\label{lem:braidedcentrallinearity}
Let us write $\mathbb{K}:=\Omega\mathcal{B}=End_{\mathcal{B}}(I)$, a finite field extension of $\mathds{k}$.
Then, $\mathcal{B}$ is linear over $\mathbb{K}$ as a (braided) fusion 1-category.
\end{Lemma}
\begin{proof}
This can be checked directly, but also follows from \cite[Remark 5.5]{KZ}.
\end{proof}

The most elementary examples of braided fusion 1-categories are given by the Drinfeld centers.
Succinctly, given a (multi)fusion 1-category $\mathcal{C}$ over $\mathds{k}$, its Drinfeld center, denoted by $\mathcal{Z}(\mathcal{C})$, is the braided monoidal 1-category whose objects are pairs $(Z,\gamma_Z)$ consisting of an object $Z$ of $\mathcal{C}$ together with a half braiding $\gamma_{Z,C}:Z\otimes C\cong C\otimes Z$ natural in the object $C$ of $\mathcal{C}$.
In fact, as observed in \cite[Corollary 2.6.8]{DSPS13}, this is a braided (multi)fusion 1-category.
Equivalently, the Drinfeld center may also be defined as the braided monoidal 1-category of $\mathcal{C}$-$\mathcal{C}$-bimodule endofunctors on $\mathcal{C}$.
Over an algebraically closed field, it is well-known that $\mathcal{Z}(\mathcal{C})$ is Morita equivalent to $\mathcal{C}^{\mathrm{mop}}\boxtimes_{\mathds{k}}\mathcal{C}$.
For later use, we record the generalization of this observation beyond the case of algebraically closed fields, noting that it requires an additional hypothesis.

\begin{Lemma}\label{lem:CenterMoritaarbitrary}
Let $\mathcal{C}$ be a fusion 1-category over $\mathds{k}$ with $\Omega\mathcal{C} = \mathds{k}$, then $\mathcal{C}$ witnesses a Morita equivalence between $\mathcal{C}^{\mathrm{mop}}\boxtimes_{\mathds{k}}\mathcal{C}$ and $\mathcal{Z}(\mathcal{C})$.
\end{Lemma}
\begin{proof}
We always have an equivalence $End_{\mathcal{C}^{\mathrm{mop}}\boxtimes_{\mathds{k}}\mathcal{C}}(\mathcal{C})\simeq \mathcal{Z}(\mathcal{C})$ as multifusion 1-categories over $\mathds{k}$.
In general, the $\mathcal{C}^{\mathrm{mop}}\boxtimes_{\mathds{k}}\mathcal{C}$-$\mathcal{Z}(\mathcal{C})$-bimodule $\mathcal{C}$ may not be invertible because the action of $\mathcal{C}^{\mathrm{mop}}\boxtimes_{\mathds{k}}\mathcal{C}$ may not be faithful.
But, as we have assumed that $\Omega\mathcal{C} = \mathds{k}$, it follows that $\mathcal{C}^{\mathrm{mop}}\boxtimes_{\mathds{k}}\mathcal{C}$ is a fusion (and not multifusion) 1-category.
This implies that the action of $\mathcal{C}^{\mathrm{mop}}\boxtimes_{\mathds{k}}\mathcal{C}$ on $\mathcal C$ is faithful, and hence that the $\mathcal{C}^{\mathrm{mop}}\boxtimes_{\mathds{k}}\mathcal{C}$-$\mathcal{Z}(\mathcal{C})$-bimodule $\mathcal{C}$ is invertible.
\end{proof}

\begin{Example}
Let us assume that $\mathds{k} = \mathbb{R}$, and take $\mathcal{C} = \mathbf{Vect}_{\mathbb{C}}$.
Then, $\mathbf{Vect}_{\mathbb{C}}\boxtimes_{\mathbb{R}}\mathbf{Vect}_{\mathbb{C}} = \mathbf{Vect}_{\mathbb{C}}\boxplus \mathbf{Vect}_{\mathbb{C}}$, whereas $\mathcal{Z}(\mathbf{Vect}_{\mathbb{C}}) = \mathbf{Vect}_{\mathbb{C}}$.
In particular, these multifusion 1-categories can not be Morita equivalent as the former is decomposable multifusion and the later is fusion.
\end{Example}

Given a fusion sub-1-category $\mathcal{A}$ of $\mathcal{B}$, we write $\mathcal{Z}_{(2)}(\mathcal{A}\subseteq\mathcal{B})$ for the centralizer of $\mathcal{A}$ in $\mathcal{B}$, that is the full fusion sub-1-category on those objects $B$ in $\mathcal{B}$ for which  $c_{X,B}\circ c_{B,X} = Id_{B\otimes X}$ for every $X$ in $\mathcal{A}$.
When $\mathcal{A} = \mathcal{B}$, we succinctly write $\mathcal{Z}_{(2)}(\mathcal{B}):=\mathcal{Z}_{(2)}(\mathcal{B}\subseteq\mathcal{B})$, which is a symmetric fusion 1-category called the symmetric, or M\"uger, center of $\mathcal{B}$.
Slightly more generally, given a braided tensor functor $F:\mathcal{B}_1\rightarrow \mathcal{B}_2$ between two braided fusion 1-categories over $\mathds{k}$, we write $\mathcal{Z}_{(2)}(F)$ for the centralizer of the image of $F$ in $\mathcal{B}_2$. We begin by analyzing how centralizers behave with respect to base change to the algebraic closure $\overline{\mathds{k}}$. Recall that $\overline{\mathcal{B}}$ denotes the base change of $\mathcal{B}$ along $\overline{\mathds{k}}/\mathds{k}$.

\begin{Lemma}\label{lem:centralizerbasechange}
Let $\mathcal{B}$ be a braided fusion 1-category over $\mathds{k}$, and $\mathcal{A}$ a full fusion sub-1-category of $\mathcal{B}$.
Then, $\mathcal{Z}_{(2)}(\overline{\mathcal{A}}\subseteq\overline{\mathcal{B}})$ is equivalent to the base change of $\mathcal{Z}_{(2)}(\mathcal{A}\subseteq\mathcal{B})$.
\end{Lemma}
\begin{proof}
Clearly, the base change of $\mathcal{Z}_{(2)}(\mathcal{A}\subseteq\mathcal{B})$ is contained in $\mathcal{Z}_{(2)}(\overline{\mathcal{A}}\subseteq\overline{\mathcal{B}})$.
It will therefore suffice to show that if $B$ in $\mathcal{B}$ is a simple object not in the centralizer of $\mathcal{A}$, then none of the summands of its image $F(B)$ in $\overline{\mathcal{B}}$ is in the centralizer of $\overline{\mathcal{A}}$.
Recall from the theory of forms developed in \cite{EG2} that $\overline{\mathcal{B}}$ comes equipped with an action of the Galois group $\mathrm{Gal}(\overline{\mathds{k}}/\mathds{k})$ whose equivariantization, i.e.\ homotopy fixed points, is $\mathcal{B}$.
There is a decomposition $$F(B) = \bigoplus_{i=1}^m \overline{B}^{\oplus n}_i,$$ where the $\overline{B}_i$ are distinct simple objects of $\overline{\mathcal{B}}$.
Moreover, the summands $\overline{B}^{\oplus n}_i$ of $F(B)$ must be permuted transitively under the action of $\mathrm{Gal}(\overline{\mathds{k}}/\mathds{k})$ as $B$ is simple.
The composite map of algebras $$End_{\mathcal{B}}(B\otimes A)\rightarrow End_{\overline{\mathcal{B}}}(F(B\otimes A))\cong End_{\overline{\mathcal{B}}}(\oplus_{i=1}^m\overline{B}^{\oplus n}_i\otimes F(A))$$ witnesses that the left hand-side is the fixed points under the canonical $\mathrm{Gal}(\overline{\mathds{k}}/\mathds{k})$-action on the right hand-side.
Let us consider the morphism $c_{A,B}\circ c_{B,A}$ in $End_{\mathcal{B}}(B\otimes A)$. By naturality of the braiding, its image under the above map of algebras splits as $$c_{F(A),F(B)}\circ c_{F(B),F(A)} = \oplus_{i=1}^m\big(c_{\overline{A},\overline{B}^{\oplus n}_i}\circ c_{\overline{B}^{\oplus n}_i,\overline{A}}\big).$$
In particular, by transitivity of the action of the Galois group, we find that $c_{A,B}\circ c_{B,A}$ is the identity if and only if $c_{\overline{A},\overline{B}^{\oplus n}_i}\circ c_{\overline{B}^{\oplus n}_i,\overline{A}}$ is the identity for any, and hence all, $i$.
The result follows.
\end{proof}

Over an algebraically closed field of characteristic zero, it is well-known that the Drinfeld center of a fusion 1-category is non-degenerate \cite[Corollary 8.20.13]{EGNO}.
We adopt the following naïve generalization of the notion of non-degeneracy.

\begin{Definition}
A braided fusion 1-category over $\mathds{k}$ is non-degenerate if the only simple object in its symmetric center is its monoidal unit.
\end{Definition}

\noindent Recall that the formation of the Drinfeld center commutes with base change to the algebraic closure by \cite[Lemma 5.1]{MS}.
The result below therefore follows from Lemma \ref{lem:centralizerbasechange}.

\begin{Corollary}\label{cor:centernondeg}
The Drinfeld center of any fusion 1-category over $\mathds{k}$ is non-degenerate.
\end{Corollary}

\begin{Remark}
A braided fusion 1-category $\mathcal{B}$ over $\mathds{k}$ is non-degenerate if its symmetric center is $\mathbf{Vect}_{\mathbb{K}}$ for some (necessarily finite) field extension $\mathbb{K}$ of $\mathds{k}$.
For the purposes of this section, this is the appropriate notion of non-degeneracy because it allows us to state that the Drinfeld center of an arbitrary fusion 1-category over $\mathds{k}$ is non-degenerate.
However, we do wish to point out that, if one wishes to think of non-degeneracy as a higher categorical invertibility condition as in \cite{BJSS}, then one would want to require the stronger condition that the symmetric center is $\mathbf{Vect}_{\mathds{k}}$.
This stronger condition will appear in section \S\ref{sec:applications} below.
\end{Remark}

Finally, we generalize \cite[Theorem 3.10]{DGNO} to arbitrary fields, but before doing so, we record the following straightforward result.

\begin{Lemma}\label{lem:=Upstairs => =Downstairs}
Let $\mathcal A$ and $\mathcal B$ be full semisimple sub-1-categories of a finite semisimple $\mathds k$-linear 1-category $\mathcal C$.
If we have $\mathcal A_{\mathbb K}\simeq\mathcal B_{\mathbb K}$ as full sub-1-categories of $\mathcal C_{\mathbb K}$ for a field extension $\mathbb K/\mathds k$, then $\mathcal A\simeq\mathcal B$.
\end{Lemma}

\begin{Lemma}\label{lem:doublecentralizer}
Let $\mathcal{B}$ be a non-degenerate braided fusion 1-category over $\mathds{k}$, and $\mathcal{A}$ a full fusion sub-1-category of $\mathcal{B}$.
Then, we have $$\mathcal{Z}_{(2)}(\mathcal{Z}_{(2)}(\mathcal{A}\subseteq \mathcal{B})\subseteq \mathcal{B})\simeq\mathcal{A}\,,$$ that is, the double centralizer of $\mathcal{A}$ in $\mathcal{B}$ is equal to $\mathcal{A}$.
\end{Lemma}
\begin{proof}
By Lemma \ref{lem:braidedcentrallinearity}, up to replacing $\mathds{k}$ by $\Omega\mathcal{B}$, we can assume that $\Omega\mathcal{B}=\mathds{k}$.
Let $\overline{\mathds{k}}$ be an algebraic closure of $\mathds{k}$ and recall that we write $\overline{\mathcal{B}}$ for the base change of $\mathcal{B}$ along the extension $\overline{\mathds{k}}/\mathds{k}$.
As we have assumed that $\Omega\mathcal{B}=\mathds{k}$, we have that $\overline{\mathcal{B}}$ is a fusion 1-category over $\overline{\mathds{k}}$.
We also write $\overline{\mathcal{A}}$ for the base change of $\mathcal{A}$ along $\overline{\mathds{k}}/\mathds{k}$, which agrees with the image of $\mathcal{A}$ along the canonical braided surjective monoidal functor $F:\mathcal{B}\rightarrow\overline{\mathcal{B}}$.
But, Lemma \ref{lem:centralizerbasechange} above implies that $\mathcal{Z}_{(2)}(\overline{\mathcal{A}}\subseteq\overline{\mathcal{B}})$ is equivalent to the base change of $\mathcal{Z}_{(2)}(\mathcal{A}\subseteq\mathcal{B})$.
The result then follows from the analogous statement over an algebraically closed field of characteristic zero given in \cite[Theorem 3.10]{DGNO} by applying Lemma \ref{lem:=Upstairs => =Downstairs}.
\end{proof}

\subsection{Technical Properties of Centers and Centralizers}

We prove various additional technical results about Drinfeld centers and centralizers that will be needed in the subsequent sections.

The next technical result is a generalization of \cite[Corollary 3.24]{DMNO}.

\begin{Proposition}\label{prop:surjectivebraided}
Let $F:\mathcal{B}_1\rightarrow \mathcal{B}_2$ be a surjective braided tensor functor between braided fusion 1-categories over $\mathds{k}$.
Then, there is a simple commutative separable algebra $A$ in $\mathcal{Z}_{(2)}(\mathcal{B}_1)$ such that $F$ is equivalent to the canonical surjective braided tensor functor $\mathcal{B}_1\rightarrow\mathbf{Mod}_{\mathcal{B}_1}(A)$.
\end{Proposition}
\begin{proof}
Let $F^*:\mathcal{B}_2\rightarrow \mathcal{B}_1$ denote the right adjoint to $F$.
Observe that $F^*$ has a canonical lax braided monoidal structure.
In particular, $A:=F^*(I)$ is a commutative algebra in $\mathcal{B}_1$ such that we have $\mathbf{Mod}_{\mathcal{B}_1}(A)\simeq \mathcal{B}_2$ as fusion 1-categories.
Furthermore, $A$ is separable by \cite[Corollary 2.6.9]{DSPS13}, and simple.
Then, it follows from the argument used in the proof of \cite[Proposition 3.22]{DMNO} that $A$ lies in $\mathcal{Z}_{(2)}(\mathcal{B}_1)$.
This yields the result.
\end{proof}

We will also use the following result, which generalizes a statement of \cite[Subsection 4.1]{DNO}.

\begin{Lemma}\label{lem:centralizersymmetriccenter}
Let $\mathcal{B}$ be a braided fusion 1-category over $\mathds{k}$.
Then, we have $$\mathcal{Z}_{(2)}(\mathcal{B})= \mathcal{Z}_{(2)}(\mathcal{B}\boxtimes_{\mathds{k}}\mathcal{B}^{\mathrm{rev}}\rightarrow\mathcal{Z}(\mathcal{B})).$$
\end{Lemma}
\begin{proof}
Note that the image of $\mathcal B\boxtimes_{\mathds k}\mathcal B^{\mathrm{rev}}\to\mathcal Z(\mathcal B)$ factors through the image of $\mathcal B\boxtimes_{\Omega\mathcal B}\mathcal B^{\mathrm{rev}}$.
Thus, without loss of generality, we may assume that $\Omega\mathcal{B} = \mathds{k}$.
Let $\overline{\mathds{k}}$ be an algebraic closure of $\mathds{k}$.
It follows from Lemma \ref{lem:centralizerbasechange} that $\mathcal{Z}_{(2)}(\overline{\mathcal{B}})$ is the base change of $\mathcal{Z}_{(2)}(\mathcal{B})$.
Moreover, by \cite[Lemma 5.1]{MS}, we have that $\mathcal{Z}(\overline{\mathcal{B}})$ agrees with the base change of $\mathcal{Z}(\mathcal{B})$.
Now, it follows from \cite[Proposition 4.3 \& Corollary 4.4]{DNO} that 
$$\mathcal{Z}_{(2)}(\overline{\mathcal{B}}\boxtimes_{\overline{\mathds{k}}}\overline{\mathcal{B}}^{\mathrm{rev}}\rightarrow \mathcal{Z}(\overline{\mathcal{B}}))=\mathcal{Z}_{(2)}(\overline{\mathcal{B}}).$$ But, $\mathcal{Z}_{(2)}(\overline{\mathcal{B}}\boxtimes_{\overline{\mathds{k}}}\overline{\mathcal{B}}^{\mathrm{rev}}\rightarrow \mathcal{Z}(\overline{\mathcal{B}}))$ coincides with the base change of $\mathcal{Z}_{(2)}(\mathcal{B}\boxtimes_{\mathds{k}}\mathcal{B}^{\mathrm{rev}}\rightarrow\mathcal{Z}(\mathcal{B}))$ thanks to Lemma \ref{lem:centralizerbasechange}.
By Lemma \ref{lem:=Upstairs => =Downstairs}, this shows that $\mathcal{Z}_{(2)}(\mathcal{B}\boxtimes_{\mathds{k}}\mathcal{B}^{\mathrm{rev}}\rightarrow\mathcal{Z}(\mathcal{B}))$ and $\mathcal{Z}_{(2)}(\mathcal{B})$ coincide as claimed.
\end{proof}

We will also use the following related observation.

\begin{Lemma}\label{lem:symcentercentralizer}
Let $\mathcal{B}$ be a non-degenerate braided fusion 1-category over $\mathds{k}$, and let $\mathcal{A}\subseteq \mathcal{B}$ be a full fusion sub-1-category. Then, we have $$\mathcal{Z}_{(2)}(\mathcal{A}) \simeq \mathcal{Z}_{(2)}(\mathcal{Z}_{(2)}(\mathcal{A}\subseteq\mathcal{B})).$$
\end{Lemma}
\begin{proof}
We have \begin{align*}\mathcal{Z}_{(2)}(\mathcal{A})&=\mathcal{A}\cap\mathcal{Z}_{(2)}(\mathcal{A}\subseteq\mathcal{B})\simeq\\ &\simeq \mathcal{Z}_{(2)}(\mathcal{Z}_{(2)}(\mathcal{A}\subseteq\mathcal{B})\subseteq\mathcal{B})\cap\mathcal{Z}_{(2)}(\mathcal{A}\subseteq\mathcal{B})=\mathcal{Z}_{(2)}(\mathcal{Z}_{(2)}(\mathcal{A}\subseteq\mathcal B)),\end{align*}
where the first and last equality follow from the definitions, and the middle one is given by Lemma \ref{lem:doublecentralizer}.
\end{proof}

The next result generalizes \cite[Theorem 3.13]{DGNO}.

\begin{Lemma}\label{lem:nondegeneracycentralizer}
Let $\mathcal{A}\hookrightarrow \mathcal{B}$ be a fully faithful inclusion of non-degenerate braided fusion 1-categories over $\mathds{k}$. There is an equivalence $$\mathcal{B}\simeq \mathcal{A}\boxtimes_{\mathbb{K}}\mathcal{Z}_{(2)}(\mathcal{A}\subseteq \mathcal{B})$$ of braided fusion 1-categories over $\mathbb{K}:=\Omega\mathcal{B}$.
\end{Lemma}
\begin{proof}
Without loss of generality, we may assume that $\mathbb{K}=\mathds{k}$. Clearly, there is a braided tensor functor $$\mathcal{A}\boxtimes_{\mathbb{K}}\mathcal{Z}_{(2)}(\mathcal{A}\subseteq \mathcal{B})\rightarrow \mathcal{B}.$$ But, this functor becomes an equivalence upon taking the base change to the algebraic closure by \cite[Theorem 3.13]{DGNO}. It follows that it must be already be an equivalence.
\end{proof}

\begin{Corollary}\label{cor:cenralizernondegenerateincenter}
Let $\mathcal{B}$ be a non-degenerate braided fusion 1-category over $\mathds{k}$. We have the following equivalence of braided fusion 1-categories over $\mathbb{K}=\Omega\mathcal{B}$ $$\mathcal{Z}(\mathcal{B})\simeq\mathcal{B}\boxtimes_{\mathbb{K}} \mathcal{B}^{\mathrm{rev}}.$$ Said differently, we have $\mathcal{Z}_{(2)}(\mathcal{B}\hookrightarrow\mathcal{Z}(\mathcal{B}))\simeq \mathcal{B}^{\mathrm{rev}}$.
\end{Corollary}
\begin{proof}
Without loss of generality, we may assume that $\mathbb{K} = \mathds{k}$. By inspection, we have that $\mathcal{B}^{\mathrm{rev}}\hookrightarrow \mathcal{Z}_{(2)}(\mathcal{B}\hookrightarrow\mathcal{Z}(\mathcal{B}))$. But, it follows from \cite[Theorem 4.9]{Sa} that $\mathrm{FPdim}(\mathcal{Z}(\mathcal{B})) = \mathrm{FPdim}(\mathcal{B})^2$ as $\mathcal{B}$ is braided and every object is therefore Galois trivial. The preceding lemma then yields $\mathrm{FPdim}(\mathcal{Z}_{(2)}(\mathcal{B}\hookrightarrow\mathcal{Z}(\mathcal{B})))=\mathrm{FPdim}(\mathcal{B})$, from which the result follows.
\end{proof}

\subsection{Symmetric Fusion 1-Categories}

We now review the theory of symmetric fusion 1-categories over $\mathds{k}$.
We will shortly state explicitly a corollary of the main result of \cite{De}. But, before doing so, we establish an auxiliary result.

\begin{Lemma}\label{lem:equivalencefinitebasechange}
Let $\mathcal{A}$ and $\mathcal{B}$ be two (braided) multifusion 1-categories over $\mathds{k}$, and let $\overline{F}:\overline{\mathcal{A}}\rightarrow \overline{\mathcal{B}}$ be a (braided) tensor functor over $\overline{\mathds{k}}$. Then, for any sufficiently large finite extension $\mathbb{L}/\mathds{k}$, there exists a (braided) tensor functor $F_{\mathbb{L}}:\mathcal{A}_{\mathbb{L}}\rightarrow \mathcal{B}_{\mathbb{L}}$ over $\mathbb{L}$ whose base change is $\overline{F}$. In particular, if we have $\overline{\mathcal{A}}\simeq \overline{\mathcal{B}}$ over $\overline{\mathds{k}}$, then $\mathcal{A}_{\mathbb{L}}\simeq \mathcal{B}_{\mathbb{L}}$ over a sufficiently large finite extension $\mathbb{L}/\mathds{k}$.
\end{Lemma}
\begin{proof}
Let $\mathbb{L}/\mathds{k}$ be a finite field extension.
We will assume that this extension is large enough so that both, for every simple object $A$ of $\mathcal{A}_{\mathbb{L}}$, we have $End_{\mathcal{A}_{\mathbb{L}}}(A)\cong\mathbb{L}$ and, for every simple object $B$ of $\mathcal{B}_{\mathbb{L}}$, we have $End_{\mathcal{B}_{\mathbb{L}}}(B)\cong\mathbb{L}$.
Said differently, the field extension $\mathbb{L}/\mathds{k}$ splits all the separable algebras of endomorphisms of the objects of $\mathcal{A}$ and $\mathcal{B}$.
In particular, it then follows that $\overline{\mathcal A} = \mathcal{A}_{\mathbb{L}}\otimes_{\mathbb{L}}\overline{\mathds{k}}$ and $\overline{\mathcal{B}} = \mathcal{B}_{\mathbb{L}}\otimes_{\mathbb{L}}\overline{\mathds{k}}$, i.e.\ it is not necessary to Cauchy complete the right hand-sides.
But then, $\mathcal{A}_{\mathbb{L}}$ and $\mathcal{B}_{\mathbb{L}}$ are non-full $\mathbb{L}$-linear (braided) multifusion sub-1-categories of $\overline{\mathcal{A}}$ and $\overline{\mathcal{B}}$. In particular, there exists an $\mathbb{L}$-linear functor $F_{\mathbb{L}}:\mathcal{A}_{\mathbb{L}}\rightarrow\mathcal{B}_{\mathbb{L}}$ such that $F_{\mathbb{L}}\otimes_{\mathbb{L}}\overline{\mathds{k}}\cong\overline{F}$. But, the coherence data for the tensor functor $\overline{F}$ is completely determined by finitely many isomorphisms in $\overline{\mathcal{B}}$, which give finitely many elements in $\overline{\mathds{k}}$. By enlarging $\mathbb{L}$ if necessary, so that it includes all of these scalars, we find that $F_{\mathbb{L}}$ admits a suitable (braided) tensor structure. This concludes the proof of the first part of the statement.
The second part follows from the observation that if $\overline{F}$ is an equivalence, then so must be $F_{\mathbb{L}}$.
\end{proof}

\begin{Theorem}[\cite{De}]\label{thm:Deligne}
Let $\mathcal{E}$ be a symmetric fusion 1-category over $\mathds{k}$.
Then, there exists a finite  field extension $\mathbb{L}$ of $\mathds{k}$ and a symmetric monoidal functor $\mathcal{E}\rightarrow \mathbf{SVect}_{\mathbb{L}}$.
\end{Theorem}
\begin{proof}
Without loss of generality, let us assume that $\Omega\mathcal{E} = \mathds{k}$.
Let us write $\overline{\mathds{k}}$ for the algebraic closure of $\mathds{k}$, and $\overline{\mathcal E}$ for the base change of $\mathcal{E}$ to $\overline{\mathds{k}}$.
It was proven in \cite{De} that there exists a fiber functor $\overline{\mathcal{E}}\rightarrow\mathbf{SVect}_{\overline{\mathds{k}}}$.
The result is therefore an immediate consequence of Lemma \ref{lem:equivalencefinitebasechange} above.
\end{proof}

\begin{Example}
Let $\mathds{k}=\mathbb{R}$, the field of real numbers.
There is a symmetric fusion 1-category $\mathbf{SVect}_{\mathbb{H}}$ of quaternionic super-vector spaces \cite[Sections 6.3 \& 6.4]{GJF0}.
Succinctly, its objects are direct sums of a vector space over $\mathbb{R}$ in even degree and a vector space over the quaternions $\mathbb{H}$ in odd degree.
Its symmetric structure is provided by the usual Koszul sign rule.
There is no monoidal functor $\mathbf{SVect}_{\mathbb{H}}\rightarrow \mathbf{SVect}_{\mathbb{R}}$.
On the other hand, there is a symmetric monoidal functor $\mathbf{SVect}_{\mathbb{H}}\rightarrow \mathbf{SVect}_{\mathbb{C}}$.
\end{Example}

Note that the (full) image of the functor $\mathcal{E}\rightarrow \mathbf{SVect}_{\mathbb{L}}$ may be either $\mathbf{Vect}_{\mathbb{L}}$ or $\mathbf{SVect}_{\mathbb{L}}$.
If the image is $\mathbf{Vect}_{\mathbb{L}}$, the symmetric fusion 1-category $\mathcal{E}$ is called Tannakian, otherwise it is called (strictly) super-Tannakian.
Now, it was argued in \cite{Del} that, over an algebraically closed field $\overline{\mathds{k}}$, every Tannakian symmetric fusion 1-category is equivalent to $\mathbf{Rep}_{\overline{\mathds{k}}}(G)$, the 1-category of finite representation of the finite group $G$.
Further, it was shown therein that every strictly super-Tannakian symmetric fusion 1-category is equivalent to $\mathbf{Rep}_{\overline{\mathds{k}}}(G,z)$, the 1-category of finite super-representation of the finite super-group $(G,z)$ with $z\neq 1$ on which $z$ acts as the parity operator.
(We note that allowing $z=1$ recovers the Tannakian case.)
For use below, we record a variant of these facts for symmetric fusion 1-categories over an arbitrary field.

\begin{Lemma}\label{lem:finitedescent}
Let $\mathcal{E}$ be a symmetric fusion 1-category over $\mathds{k}$ such that $\Omega\mathcal{E} = \mathds{k}$.
There exists a finite normal field extension $\mathbb{L}$ such that the base change $\mathcal{E}_{\mathbb{L}}$ is equivalent to either $\mathbf{Rep}_{\mathbb{L}}(G)$ for some finite group $G$ if $\mathcal{E}$ is Tannakian or $\mathbf{Rep}_{\mathbb{L}}(G,z)$ for a finite super-group $(G,z)$ if $\mathcal{E}$ is super-Tannakian.
\end{Lemma}
\begin{proof}
By the preceding discussion, this follows from \cite{Del} and Lemma \ref{lem:equivalencefinitebasechange}.
\end{proof}

\subsection{Picard Groups}

Over an algebraically closed field of characteristic zero, the Picard group of a braided fusion 1-category was introduced in \cite{ENO2}.
This construction can be generalized to any braided fusion 1-category $\mathcal{B}$ over an arbitrary field $\mathds{k}$.
More precisely, we write $$\mathcal{P}ic(\mathcal{B}) := \mathbf{Mod}(\mathcal{B})^{\times}$$ for the space of invertible $\mathcal{B}$-module 1-categories and their morphisms.
The relative Deligne tensor product endows this space with the structure of a group-like topological monoid.
We use $\mathrm{Pic}(\mathcal{B}) := \pi_0(\mathcal{P}ic(\mathcal{B}))$ to denote the group of (equivalence classes of) invertible $\mathcal{B}$-module 1-categories.
Alternatively, $\mathcal{P}ic(\mathcal{B})$ corresponds to the space of Azumaya algebras in $\mathcal{B}$, invertible bimodules between them, and bimodule morphisms.
In particular, as remarked in \cite[Proposition 4.9]{ENO2}, we have $$\mathrm{Pic}(\mathbf{Vect}_{\mathds{k}}) = \mathrm{Br}(\mathds{k}),$$ the Picard group of $\mathbf{Vect}_{\mathds{k}}$ coincides with the Brauer group of $\mathds{k}$.
But, it is a classical result that the Brauer group of a field is always torsion.
We generalize this fact to the Picard group of an arbitrary symmetric fusion 1-category over $\mathds{k}$.

\begin{Proposition}\label{prop:Pictorsion}
Let $\mathcal{E}$ be a symmetric fusion 1-category over $\mathds{k}$.
The Picard group $\mathrm{Pic}(\mathcal{E})$ is torsion.
\end{Proposition}
\begin{proof}
Without loss of generality, we may assume that $\Omega\mathcal{E} = \mathds{k}$.
Namely, the Picard group does not depend on the ground field.
Then, it follows from Lemma \ref{lem:finitedescent} that there exists a finite normal field extension $\mathbb{K}$ such that the base change $\mathcal{F}:=\mathcal{E}_{\mathbb{K}}$ is equivalent to the 1-category of representations over $\mathbb{K}$ of a finite (super-)group.
Moreover, writing $\Gamma := \mathrm{Gal}(\mathbb{K}/\mathds{k})$, a finite group by construction, there is an action of $\Gamma$ on $\mathcal{F}$ by symmetric monoidal $\mathds{k}$-linear autoequivalences such that $\mathcal{F}^{\Gamma} \simeq \mathcal{E}$.

Firstly, we observe that the abelian group $\mathrm{Pic}(\mathcal{F})$ is torsion.
On the one hand, if $\mathcal{F}$ is Tannakian, there exists a finite group $G$ and an equivalence $\mathcal{F}\simeq \mathbf{Rep}_{\mathbb{K}}(G)$ by choice of $\mathbb{K}$.
But, it was shown in \cite[Theorem 1.12]{Long} that $\mathrm{Pic}(\mathbf{Rep}_{\mathbb{K}}(G))=\mathrm{Br}(\mathbb{K})\oplus H^2_{Sw}(\mathbb{K}[G];\mathbb{K})$, where $H^2_{Sw}(\mathbb{K}[G];\mathbb{K})$ denotes Sweedler's cohomology for the Hopf algebra $\mathbb{K}[G]$.
Moreover, in this case, we have $H^2_{Sw}(\mathbb{K}[G];\mathbb{K})\cong H^2(G;\mathbb{K}^\times)$ by \cite[Theorem 3.1]{Swe}.
It follows that $\mathrm{Pic}(\mathbf{Rep}_{\mathbb{K}}(G))$ is torsion.
On the other hand, if $\mathcal{F}$ is strictly super-Tannakian, then there exists a finite super-group $(G,z)$ with $z\neq 1$ and an equivalence $\mathcal{F}\simeq \mathbf{Rep}_{\mathbb{K}}(G,z)$.
It follows from \cite[Theorems 2.9 \& 2.10]{Car} that $\mathrm{Pic}(\mathbf{Rep}_{\mathbb{K}}(G,z))$ is torsion.

Secondly, we consider the canonical group homomorphism $\mathrm{Pic}(\mathcal{E})\rightarrow \mathrm{Pic}(\mathcal{F})$.
As the target is torsion by the discussion in the preceding paragraph, it will be enough to show that the kernel is torsion.
Let $A$ be an Azumaya algebra in $\mathcal{E}$ whose image under $\mathcal{E}\rightarrow \mathcal{F}$ is Morita trivial.
This is equivalent to the data of a non-zero object $V$ of $\mathcal{F}$, together with a $\Gamma$-fixed point structure on the algebra $\underline{End}_{\mathcal{F}}(V)=V\otimes V^*$ in $\mathcal{F}$.
Moreover, the algebra $A$ in $\mathcal{E}\simeq \mathcal{F}^{\Gamma}$ is Morita trivial precisely if the $\Gamma$-fixed point structure on $\underline{End}_{\mathcal{F}}(V)$ can be lifted to a $\Gamma$-fixed point structure on $V$.
We will argue that there are two obstructions to this lift, which are measured by the torsion abelian groups $H^1(\Gamma; H^1(G;\mathbb{K}^\times))$ and $H^2(\Gamma;\mathbb{K}^\times)$.
This will complete the proof of the proposition.

Let us now assume that the algebra $\underline{End}_{\mathcal{F}}(V)$ in $\mathcal{F}$ comes equipped with a $\Gamma$-fixed point structure (under the action of $\Gamma$ on $\mathcal{F}$).
In particular, there are isomorphisms of algebras $$f_{\gamma}:(\gamma\cdot V)\otimes (\gamma\cdot V)^*\cong \gamma\cdot \underline{End}_{\mathcal{F}}(V)\cong \underline{End}_{\mathcal{F}}(V)\cong V\otimes V^*.$$ 
As plain fusion 1-categories, we have $\mathcal{F}\simeq \mathbf{Rep}_{\mathbb{K}}(G)$, and $\mathcal{F}$ therefore admits a canonical spherical structure such that $\mathrm{Tr}(Id_V) = \mathrm{dim}_{\mathbb{K}}(V)$, the dimension of the underlying $\mathbb{K}$-vector space.
We emphasize that, when $\mathcal{F}$ is strictly super-Tannakian, the spherical structure we are considering is not compatible with the symmetric structure.
Taking mates of the morphisms $f_{\gamma}$,
we obtain endomorphisms $p_{\gamma}:V^*\otimes(\gamma\cdot V)\rightarrow V^*\otimes (\gamma\cdot V)$ in $\mathcal{F}$.
It follows from the fact that $f_{\gamma}$ is a homomorphism of algebras that $p_{\gamma}\circ p_{\gamma} = \mathrm{dim}_{\mathbb{K}}(V)p_{\gamma}$ and $\mathrm{Tr}_{\mathbb{K}}(p_{\gamma}) = \mathrm{dim}_{\mathbb{K}}(V)$.
But, $\mathrm{dim}_{\mathbb{K}}(V)\neq 0$ as $V$ is non-zero, so that $p_{\gamma}/\mathrm{dim}_{\mathbb{K}}(V)$ is an idempotent on $V^*\otimes (\gamma\cdot V)$ of trace $1$.
This implies that $p_{\gamma}/\mathrm{dim}_{\mathbb{K}}(V)$ corresponds to the inclusion of a one dimensional $G$-module $J_{\gamma}$ over $\mathbb{K}$. 
Said differently, $J_{\gamma}$ is an invertible objects of $\mathcal{F}$, and these form the group $H^1(G;\mathbb{K}^\times)$.
The assignment $\gamma\mapsto J_{\gamma}$ therefore defines a map $\Gamma\rightarrow H^1(G;\mathbb{K}^\times)$, and one checks easily that it defines a class $\xi$ in $H^1(\Gamma; H^1(G;\mathbb{K}^\times))$, which is precisely the first obstruction group mentioned above.
Now, using the symmetric structure on $\mathcal{F}$, we can endow the tensor power $\underline{End}_{\mathcal{F}}(V)^{\otimes n}$ with an algebra structure.
It is straightforward to check that the corresponding class in $H^1(\Gamma; H^1(G;\mathbb{K}^\times))$ is $\xi^n$, and, in particular, does not actually depend on the symmetric structure of $\mathcal{F}$.
Up to replacing $\underline{End}_{\mathcal{F}}(V)$ by a tensor power, we may therefore assume that the obstruction $\xi$ vanishes.

Finally, if $J_{\gamma} = I$ for all $\gamma$ in $\Gamma$, then $f_{\gamma}$ must factor as $l_{\gamma}\otimes (l_{\gamma}^*)^{-1}$ for a unique isomorphism $l_{\gamma}:\gamma\cdot V\cong V$ in $\mathcal{F}$.
The isomorphisms $l_{\gamma}$ assemble into a $\Gamma$-fixed point structure on $V$ if and only if the squares in $\mathcal{F}$ $$\begin{tikzcd}
\delta\cdot (\gamma\cdot V) \arrow[d, "\rotatebox{90}{$\cong$}"'] \arrow[r, "\delta\cdot l_{\gamma}"] & \gamma\cdot V \arrow[d, "l_{\gamma}"] \\
(\delta\gamma)\cdot V \arrow[r, "l_{\delta\gamma}"']                                & V                               \end{tikzcd}$$ commute for all $\gamma$ and $\delta$.
But, by hypothesis, the squares $$\begin{tikzcd}
\big(\delta\cdot (\gamma\cdot V)\big)\otimes\big(\delta\cdot (\gamma\cdot V)\big)^* \arrow[d, "\rotatebox{90}{$\cong$}"'] \arrow[r, "\delta\cdot f_{\gamma}"] & (\gamma\cdot V)\otimes (\gamma\cdot V)^* \arrow[d, "f_{\gamma}"] \\
\big((\delta\gamma)\cdot V\big)\otimes \big((\delta\gamma)\cdot V\big)^* \arrow[r, "f_{\delta\gamma}"']                                     & V\otimes V^*
\end{tikzcd}$$ of algebra isomorphisms in $\mathcal{F}$ do commute, so that the failure of commutativity of the first square above for fixed $\delta$ and $\gamma$ is given by a scalar $\lambda_{\delta,\gamma}$ in $\mathbb{K}^{\times}$.
Further, the assignment $(\delta,\gamma)\mapsto\lambda_{\delta,\gamma}$ defines an element of $H^2(\Gamma;\mathbb{K}^{\times})$, the second obstruction group mentioned above.
This completes the proof of the proposition.
\end{proof}

We now derive a non-commutative version of the last proposition above.
More precisely, the Picard group of an arbitrary braided fusion 1-category over an algebraically closed field is in general not abelian, but it is always finite by \cite[Theorem 4.15]{ENO2}.
Further, while the Picard group of a braided fusion 1-category (over an arbitrary field) is in general not finite, it is still natural to ask whether it is nevertheless torsion.
In fact, it will be important for our subsequent applications to prove that such Picard groups satisfy a stronger property:\ Namely, observe that every torsion abelian group is indfinite, also referred to as locally finite, meaning that every finitely generated subgroup is finite.
We will show shortly that the Picard group of an arbitrary braided fusion 1-category satisfies this property, but, in order to do so, we will appeal to the following technical lemma.

\begin{Lemma}\label{lem:indfiniteextension}
Let $$1\rightarrow K\rightarrow G\xrightarrow{\phi} Q\rightarrow 1$$ be an exact sequence of groups, that is, $K$ is the kernel of $\phi$.
If $K$ and $Q$ are indfinite, then so is $G$.
\end{Lemma}
\begin{proof}
Recall from \cite{EML} that such an arbitrary extension of groups is classified by the data of two maps of sets $\psi:Q\rightarrow \mathrm{Aut}(K)$ and $\chi:Q\times Q\rightarrow K$ satisfying various conditions.
Moreover, the group $G$ may then be identified with the direct product $Q\times K$ equipped with the group structure $$(q_1,k_1)\cdot (q_2,k_2) = (q_1q_1,\chi(q_1,q_1)\big(\psi(q_2)^{-1}(k_1)\big)k_2).$$ We will use this identification liberally in the remainder of the proof.

It is enough to show that every finitely generated subgroup of $G$ is finite.
Let $(q_i,k_i)$ with $i=1,...,n$ be a finite collection of elements of $G$.
As $Q$ is indfinite by assumption, we find that $q_i$ generate a finite subgroup of $Q$ that we denote by $P$.
Moreover, let us consider the subgroup $J$ of $K$ generated by the finitely many elements $k_i$, $\psi(p_1)^{-1}(k_i)$, and $\chi(p_1,p_2)$ where $i=1,...,n$, and $p_1,p_2$ range over all pairs of elements of $P$.
As $K$ is indfinite by hypothesis, the subgroup $J$ must be finite.
Clearly, $\chi(p_1,p_2)$ lies in $J$ for every $p_1,p_2$ in $P$.
Moreover, it follows from the equations satisfied by $\psi$ and $\chi$ that $\psi(p)$ restricts to an automorphism of $J$ for every $p$ in $P$.
Thence, we find that the subset $P\times J$ of $Q\times K$ is closed under the group structure.
But, $P\times J$ is finite as both $P$ and $J$ are, and contains the elements $(q_i,k_i)$.
This concludes the proof.
\end{proof}

\begin{Theorem}\label{thm:Picardgrouptorsion}
Let $\mathcal{B}$ be a braided fusion 1-category over $\mathds{k}$.
The Picard group $\mathrm{Pic}(\mathcal{B})$ is indfinite, i.e.\ every finitely generated subgroup is finite.
\end{Theorem}
\begin{proof}
Firstly, note that the group $\mathrm{Pic}(\mathcal{B})$ does not depend on $\mathds{k}$.
Namely, $\mathrm{Pic}(\mathcal{B})$ is the group of Morita equivalence classes of Azumaya algebras in $\mathcal{B}$, so that we may assume that $\Omega\mathcal{B} = \mathds{k}$.
Secondly, it follows from the generalization of \cite[Remark 5.5]{DN} to arbitrary fields of characteristic zero, that there is a fiber sequence $$\mathcal{P}ic^{\mathrm{br}}(\mathcal{B})\rightarrow\mathcal{P}ic(\mathcal{B})\rightarrow \mathcal{A}ut^{\mathrm{br}}_{\mathds{k}}(\mathcal{B})$$ of group-like topological monoids, where $\mathcal{A}ut^{\mathrm{br}}_{\mathds{k}}(\mathcal{B})$ denotes the space of $\mathds{k}$-linear braided tensor autoequivalences of $\mathcal{B}$.
Associated to the fiber sequence above, there is a long exact sequence of homotopy groups \begin{align}\notag 0 &\rightarrow \mathrm{Inv}(\mathcal{E})\rightarrow \mathrm{Inv}(\mathcal{B})\rightarrow \mathrm{Aut}^{\mathrm{br}}_{\mathds{k}}(Id_{\mathcal{B}})\rightarrow\\ &\rightarrow \mathrm{Pic}^{\mathrm{br}}(\mathcal{B})\rightarrow \mathrm{Pic}(\mathcal{B})\rightarrow \mathrm{Aut}^{\mathrm{br}}_{\mathds{k}}(\mathcal{B}).\label{eq:LESPicard}\end{align}
By Lemma \ref{lem:indfiniteextension}, in order to show that $\mathrm{Pic}(\mathcal{B})$ is indfinite, it will be enough to show that both $\mathrm{Pic}^{\mathrm{br}}(\mathcal{B})$ and $\mathrm{Aut}^{\mathrm{br}}_{\mathds{k}}(\mathcal{B})$ are indfinite.

Thirdly, we argue that $\mathrm{Pic}^{\mathrm{br}}(\mathcal{B})$ is indfinite.
Let us write $\mathcal{E}:=\mathcal{Z}_{(2)}(\mathcal{B})$, a symmetric fusion 1-category over $\mathds{k}$.
The group $\mathrm{Pic}(\mathcal{E})$ is torsion by Proposition \ref{prop:Pictorsion}, and therefore indfinite as it is also abelian.
Now, recall from \cite[Theorem 4.11]{DN} that $\mathcal{P}ic^{\mathrm{br}}(\mathcal{B})\simeq \mathscr{Z}(\mathbf{Mod}(\mathcal{B}))^{\times}$.
This shows that the group $\mathrm{Pic}^{\mathrm{br}}(\mathcal{B})$ is abelian.
Then, as $\Omega\mathscr{Z}(\mathbf{Mod}(\mathcal{B}))\simeq \mathcal{E}$, the inclusion of the connected component of the identity $\mathbf{Mod}(\mathcal{E})\hookrightarrow \mathscr{Z}(\mathbf{Mod}(\mathcal{B}))$ induces an inclusion of abelian group $\mathrm{Pic}(\mathcal{E})\rightarrow \mathrm{Pic}^{\mathrm{br}}(\mathcal{B})$.
Its cokernel is the group of connected components of $\mathscr{Z}(\mathbf{Mod}(\mathcal{B}))$ containing an invertible object.
But, this group is necessarily finite as $\mathscr{Z}(\mathbf{Mod}(\mathcal{B}))$ has finitely many connected components.
It therefore follows that $\mathrm{Pic}^{\mathrm{br}}(\mathcal{B})$ is both torsion and abelian, so that it is indeed indfinite.

Fourthly, it only remains to prove that the group $\mathrm{Aut}^{\mathrm{br}}_{\mathds{k}}(\mathcal{B})$ is indfinite.
To this end, recall that the group $\mathrm{Aut}^{\mathrm{br}}_{\overline{\mathds{k}}}(\overline{\mathcal{B}})$ of braided tensor autoequivalences of $\overline{\mathcal{B}}$, the base change of $\mathcal{B}$ to the algebraic closure $\overline{\mathds{k}}$ of $\mathds{k}$, is finite (see for instance \cite[Theorem 9.1.5]{EGNO}).
Thanks to Lemma \ref{lem:indfiniteextension}, it is therefore enough to show that the kernel of the group homomorphism $\mathrm{Aut}^{\mathrm{br}}_{\mathds{k}}(\mathcal{B})\rightarrow \mathrm{Aut}^{\mathrm{br}}_{\overline{\mathds{k}}}(\overline{\mathcal{B}})$ is indfinite.
Now, given any braided tensor autoequivalence $F:\mathcal{B}\simeq \mathcal{B}$ whose base change $\overline{F}$ to $\overline{\mathds{k}}$ is equivalent to the identity, there exists a finite field extension $\mathbb{K}$ of $\mathds{k}$ such that the base change $F_{\mathbb{K}}$ is equivalent to the identity on $\mathcal{B}_{\mathbb{K}}$.
Namely, a monoidal natural equivalence only involves the data of finitely many scalars in $\overline{\mathds{k}}$.

Finally, the result will follow provided that we can show that the kernel of the map $\mathrm{Aut}^{\mathrm{br}}_{\mathds{k}}(\mathcal{B})\rightarrow \mathrm{Aut}^{\mathrm{br}}_{\mathbb{K}}(\mathcal{B}_{\mathbb{K}})$ is indfinite for every finite field extension $\mathbb{K}$ of $\mathds{k}$.
This follows from the theory of forms \cite{EG2}.
More precisely, we may without loss of generality assume that $\mathbb{K}$ is a finite Galois extension of $\mathds{k}$.
Then, the Galois group $\mathrm{Gal}(\mathbb{K}/\mathds{k})$, which is finite, canonically acts on $\mathcal{B}_{\mathbb{K}}$.
The homotopy fixed points under this action are given by $\mathcal{B}$.
Furthermore, the data of a braided monoidal autoequivalence $F:\mathcal{B}\simeq \mathcal{B}$ whose base change $F_{\mathbb{K}}$ is the identity $Id_{\mathcal{B}_{\mathbb{K}}}$ is precisely the data of a form on $Id_{\mathcal{B}_{\mathbb{K}}}$.
But, upon inspecting the definitions, we find that a form on $Id_{\mathcal{B}_{\mathbb{K}}}$ is equivalent to the data of a group homomorphism $\mathrm{Gal}(\mathbb{K}/\mathds{k})\rightarrow \mathrm{Aut}^{\mathrm{br}}_{\mathds{k}}(Id_{\mathcal{B}_{\mathbb{K}}})$.
It is therefore sufficient to show that the abelian group $\mathrm{Aut}^{\mathrm{br}}_{\mathds{k}}(Id_{\mathcal{B}_{\mathbb{K}}})=\mathrm{Aut}^{\mathrm{br}}_{\mathbb{K}}(Id_{\mathcal{B}_{\mathbb{K}}})$ is indfinite, or, equivalently, torsion.
This follows from the long exact sequence \eqref{eq:LESPicard} for $\mathcal{B}_{\mathbb{K}}$ as $\mathrm{Pic}^{\mathrm{br}}(\mathcal{B}_{\mathbb{K}})$ is abelian and torsion by the discussion in the first paragraph above.
This concludes the proof of the proposition.
\end{proof}

\begin{Remark}
We suspect that Theorem \ref{thm:Picardgrouptorsion} can be generalized further.
More precisely, we expect that the Picard group of any finite braided tensor 1-category over an arbitrary field (of arbitrary characteristic) is always indfinite.
\end{Remark}

\section{The Morita Theory of Connected Compact Semisimple Tensor 2-Categories}\label{sec:connected}

Working over an arbitrary field $\mathds{k}$ of characteristic zero, we generalize the results of \cite[Section 3]{D9} to the non-algebraically closed case.
More precisely, we generalize the concept of Witt equivalence between braided fusion 1-categories over an algebraically closed field introduced in \cite{DMNO,DNO} to braided fusion 1-category over an arbitrary field $\mathds{k}$.
We then compare this notion with that of Morita equivalence between connected compact semisimple tensor 2-categories.

\subsection{Witt Equivalence}

Let $\mathcal{B}$ be a braided fusion 1-category over $\mathds{k}$.
Let $\mathcal{C}$ be a $\mathcal{B}$-central fusion 1-category over $\mathds{k}$, with braided tensor functor $F:\mathcal{B}\rightarrow\mathcal{Z}(\mathcal{C})$.
We write $\mathcal{Z}(\mathcal{C},\mathcal{B})$ for the centralizer of the image of $F$ in $\mathcal{Z}(\mathcal{C})$, that is, we set $$\mathcal{Z}(\mathcal{C},\mathcal{B}):= \mathcal{Z}_{(2)}(F:\mathcal{B}\rightarrow \mathcal{Z}(\mathcal{C})).$$ Explicitly, $\mathcal{Z}(\mathcal{C},\mathcal{B})$ is the full fusion sub-1-category of $\mathcal{Z}(\mathcal{C})$ on those object $Z$ of $\mathcal{Z}(\mathcal{C})$ such that $c_{F(X),Z}\circ c_{Z,F(X)} = Id_{Z\otimes F(X)}$ for every object $X$ in $\mathcal{B}$.

\begin{Remark}\label{rem:braidedcentrallinearity2}
Lemma \ref{lem:braidedcentrallinearity} records the fact that, with $\mathbb{K}:=\Omega\mathcal{B}=End_{\mathcal{B}}(I)$, $\mathcal{B}$ is linear over $\mathbb{K}$.
Now, if $\mathcal{C}$ is a $\mathcal{B}$-central fusion 1-category, then $\mathcal{C}$ admits a unique $\mathbb{K}$-linear structure such that $F$ is $\mathbb{K}$-linear.
Namely, $\Omega F$ is an embedding of fields and $\mathcal{C}$ is linear over $\Omega\mathcal{Z}(\mathcal{C})$ as noted in \cite[Remark 5.5]{KZ}.
In particular, the braided fusion 1-category $\mathcal{Z}(\mathcal{C},\mathcal{B})$ is also $\mathbb{K}$-linear.
\end{Remark}

\begin{Definition}\label{def:unpointedWitt}
Let $\mathcal{B}_1$ and $\mathcal{B}_2$ be two braided fusion 1-categories over $\mathds{k}$.
We say that $\mathcal{B}_1$ and $\mathcal{B}_2$ are Witt equivalent provided that there exists a fusion 1-category $\mathcal{C}$ over $\mathds{k}$,
together with a fully faithful braided embedding $\mathcal{B}_1\hookrightarrow\mathcal{Z}(\mathcal{C})$ and a braided
monoidal equivalence $\mathcal{B}_2\simeq \mathcal{Z}(\mathcal{C},\mathcal{B}_1)^{\mathrm{rev}}$ over $\mathds{k}$.
\end{Definition}

As an application of Lemma \ref{lem:symcentercentralizer}, we obtain the following useful observation.

\begin{Corollary}
Let $\mathcal{B}_1$ and $\mathcal{B}_2$ be two braided fusion 1-categories over $\mathds{k}$. If they are Witt equivalent, then $\mathcal{Z}_{(2)}(\mathcal{B}_1)\simeq \mathcal{Z}_{(2)}(\mathcal{B}_2)$ as symmetric fusion 1-categories over $\mathds{k}$.
\end{Corollary}

\begin{Remark}\label{rem:relationWittgroups}
Working over an algebraically closed field of characteristic zero $\overline{\mathds{k}}$, the concept of Witt equivalence between non-degenerate braided fusion 1-categories over $\overline{\mathds{k}}$, that is, braided fusion 1-categories with symmetric center $\mathbf{Vect}_{\overline{\mathds{k}}}$, was introduced in \cite{DMNO}.
There is a corresponding group $\mathrm{Witt} = \mathrm{Witt}(\mathbf{Vect}_{\overline{\mathds{k}}})$, called the Witt group, which is the quotient of the monoid of non-degenerate braided fusion 1-categories by the submonoid of Drinfeld centers.
Generalizing this notion, a definition of Witt equivalence for braided fusion 1-categories over $\overline{\mathds{k}}$ whose symmetric center is explicitly identified with a fixed symmetric fusion 1-category $\mathcal{E}$ was given in \cite{DNO}.
The corresponding group is denoted by $\mathrm{Witt}(\mathcal{E})$.
These constructions can readily be generalized to the case of an arbitrary field $\mathds{k}$.
That is, for any symmetric fusion 1-category $\mathcal{E}$ over $\mathds{k}$, one can consider the Witt group $\mathrm{Witt}(\mathcal{E})$ of braided fusion 1-categories with symmetric center identified with $\mathcal{E}$.
The group $\mathrm{Aut}^{\mathrm{br}}(\mathcal{E})$ of symmetric monoidal autoequivalences of $\mathcal{E}$ acts on $\mathrm{Witt}(\mathcal{E})$, and we write $\widehat{\mathrm{Witt}}(\mathcal{E})$ for the corresponding set of orbits under this action.
Generalizing \cite[Lemma 3.1.2]{D9} to an arbitrary field, two braided fusion 1-categories $\mathcal{B}_1$ and $\mathcal{B}_2$ are Witt equivalent in the sense of definition \ref{def:unpointedWitt} above if and only if they have the same symmetric center $\mathcal{E}$ and yield the same class in $\widehat{\mathrm{Witt}}(\mathcal{E})$.
\end{Remark}

\begin{Remark}
We give a higher categorical perspective on the notions of this section.
A symmetric monoidal Morita 4-category $\mathbf{BrFus}_{\mathds{k}}$ of braided fusion 1-categories over $\mathds{k}$ was constructed in \cite{BJS}.
One can show that the above notion of Witt equivalence coincides with the property that two braided fusion 1-categories $\mathcal{B}_1$ and $\mathcal{B}_2$ over $\mathds{k}$ are equivalent when viewed as objects of $\mathbf{BrFus}_{\mathds{k}}$, with equivalence supplied by the $\mathcal{B}_1\boxtimes\mathcal{B}_2^{\mathrm{rev}}$-central fusion 1-category $\mathcal{C}$.
Over an algebraically closed field $\overline{\mathds{k}}$, this was shown in \cite{BJSS} for Witt equivalence between non-degenerate braided fusion 1-categories.
In particular, they showed that the group of invertible objects of $\mathbf{BrFus}_{\overline{\mathds{k}}}$ is precisely $\mathrm{Witt}(\mathbf{Vect}_{\overline{\mathds{k}}})$.
This can be generalized to an arbitrary field.
Yet more generally, for any fixed symmetric fusion 1-category $\mathcal{E}$ over $\mathds{k}$, one can consider the symmetric monoidal Morita 4-category $\mathbf{BrFus}_{\mathcal{E}}$ of braided fusion 1-categories enriched over $\mathcal{E}$ (using \cite[Section 2]{BJSS}).
The group of (equivalence classes of) invertible objects of $\mathbf{BrFus}_{\mathcal{E}}$ is exactly $\mathrm{Witt}(\mathcal{E})$.
\end{Remark}

\subsection{Comparison with Morita Equivalence}

We now compare the notion of Witt equivalence between braided fusion 1-categories introduced above with that of Morita equivalence between the corresponding connected compact semisimple tensor 2-categories.

\begin{Theorem}\label{thm:MoritaWitt}
Let $\mathcal{B}_1$ and $\mathcal{B}_2$ be two braided fusion 1-categories over $\mathds{k}$.
Then, $\mathbf{Mod}(\mathcal{B}_1)$ and $\mathbf{Mod}(\mathcal{B}_2)$ are Morita equivalent if and only if $\mathcal{B}_1$ and $\mathcal{B}_2$ are Witt equivalent.
\end{Theorem}

\noindent In essence, the proof follows the argument used to establish \cite[Theorem 3.1.4]{D9}.
We give some details because we must appeal to generalizations of the results used therein.
In particular, we will need the following technical lemmas.

\begin{Lemma}\label{lem:loopsdualconnected}
Let $\mathcal{B}$ be a braided fusion 1-category over $\mathds{k}$.
Let $\mathcal{C}$ be a $\mathcal{B}$-central fusion 1-category over $\mathds{k}$.
Then, the Morita dual $\mathfrak{C}$ to $\mathbf{Mod}(\mathcal{B})$ with respect to $\mathbf{Mod}(\mathcal{C})$ has loops $\Omega\mathfrak{C}$ given by $\mathcal{Z}(\mathcal{C},\mathcal{B})^{\mathrm{rev}}$, the centralizer of the image of $\mathcal{B}$ in $\mathcal{Z}(\mathcal{C})$.
\end{Lemma}
\begin{proof}
This is exactly as in \cite[Lemma 3.2.1]{D9}.
\end{proof}

\begin{Corollary}\label{cor:loopscenter}
Let $\mathcal{B}$ be a braided fusion 1-category over $\mathds{k}$.
There is an equivalence of symmetric fusion 1-categories $$\Omega\mathscr{Z}(\mathbf{Mod}(\mathcal{B}))\simeq \mathcal{Z}_{(2)}(\mathcal{B}).$$
\end{Corollary}
\begin{proof}
This follows from the lemma above using that the dual to $\mathbf{Mod}(\mathcal{B})\boxtimes_{\mathds{k}} \mathbf{Mod}(\mathcal{B})^{\mathrm{mop}}$ with respect to $\mathbf{Mod}(\mathcal{B})$ is $\mathscr{Z}(\mathbf{Mod}(\mathcal{B}))$ together with the equivalence $\mathcal{Z}(\mathcal{B}, \mathcal{B}\boxtimes_{\mathds{k}}\mathcal{B}^{\mathrm{rev}})\simeq \mathcal{Z}_{(2)}(\mathcal{B})$ from Lemma \ref{lem:centralizersymmetriccenter}.
Alternatively, this result can also be proven directly by inspecting the definitions.
\end{proof}

\begin{Lemma}
Let $\mathcal{B}$ be a braided fusion 1-category over $\mathds{k}$, and $\mathcal{C}$ a $\mathcal{B}$-central fusion 1-category.
If the braided tensor functor $F:\mathcal{B}\hookrightarrow \mathcal{Z}(\mathcal{C})$ is fully faithful, then $\mathbf{Bimod}_{\mathbf{Mod}(\mathcal{B})}(\mathcal{C})$ is connected.
\end{Lemma}
\begin{proof}
By Remark \ref{rem:braidedcentrallinearity2}, we may assume that $\Omega\mathcal{B}=\mathds{k}$.
As we have assumed that $F$ is fully faithful, observe that this implies that $\Omega\mathcal{Z}(\mathcal{C})=\Omega\mathcal{B} = \mathds{k}$.
In particular, the base change $\overline{F}:\overline{\mathcal{B}}\hookrightarrow\mathcal{Z}(\overline{\mathcal{C}})$ is a fully faithful inclusion of braided fusion 1-categories over $\overline{\mathds{k}}$.
In particular, the fusion 2-category $\mathbf{Bimod}_{\mathbf{Mod}(\overline{\mathcal{B}})}(\overline{\mathcal{C}})$ is connected by \cite[Lemma 3.2.5]{D9}.
But, it follows from the proof of \cite[Theorem 3.1.4]{D10} that the base change of the compact semisimple tensor 2-category $\mathbf{Bimod}_{\mathbf{Mod}(\mathcal{B})}(\mathcal{C})$ is $\mathbf{Bimod}_{\mathbf{Mod}(\overline{\mathcal{B}})}(\overline{\mathcal{C}})$.
We therefore find that $\mathbf{Bimod}_{\mathbf{Mod}(\mathcal{B})}(\mathcal{C})$ must be connected.
\end{proof}

\begin{Corollary}\label{cor:dualconnected}
Let $\mathcal{B}$ be a braided fusion 1-category over $\mathds{k}$, and $\mathcal{C}$ a $\mathcal{B}$-central fusion 1-category over $\mathds{k}$ such that the braided tensor functor $F:\mathcal{B}\rightarrow \mathcal{Z}(\mathcal{C})$ is fully faithful.
Then, the Morita dual tensor 2-category to $\mathbf{Mod}(\mathcal{B})$ with respect to $\mathbf{Mod}(\mathcal{C})$ is $\mathbf{Mod}(\mathcal{Z}(\mathcal{C},\mathcal{B})^{\mathrm{rev}})$.
\end{Corollary}

\vspace{1mm}

\begin{proof}[Proof of Thm.
\ref{thm:MoritaWitt}]
By Corollary \ref{cor:dualconnected}, it follows that if $\mathcal{B}_1$ and $\mathcal{B}_2$ are Witt equivalent braided fusion 1-categories over $\mathds{k}$, then $\mathbf{Mod}(\mathcal{B}_1)$ and $\mathbf{Mod}(\mathcal{B}_2)$ are Morita equivalent.
This establishes one direction, it therefore only remains to consider the converse.
Assume that $\mathbf{Mod}(\mathcal{B}_1)$ and $\mathbf{Mod}(\mathcal{B}_2)$ are Morita equivalent compact semisimple tensor 2-categories.
Then, they have equivalent Drinfeld centers by Theorem \ref{thm:MoritaCenter}.
It therefore follows from Lemma \ref{cor:loopscenter} that $\mathcal{B}_1$ and $\mathcal{B}_2$ have the same symmetric center $\mathcal{E}$.
By hypothesis, there exists a $\mathcal{B}_1$-central fusion 1-category $\mathcal{C}$ together with a braided tensor functor $F:\mathcal{B}_1\rightarrow\mathcal{Z}(\mathcal{C})$ such that $\mathbf{Bimod}_{\mathbf{Mod}(\mathcal{B}_1)}(\mathcal{C})\simeq \mathbf{Mod}(\mathcal{B}_2)$ as compact semisimple tensor 2-categories.
Let us write $\mathcal{A}$ for the (full) image of $\mathcal{B}_1$ in $\mathcal{Z}(\mathcal{C})$.
It follows from Lemma \ref{lem:loopsdualconnected} that $$\mathcal{Z}(\mathcal{C},\mathcal{B}_1)=\mathcal{Z}(\mathcal{C},\mathcal{A})\simeq \mathcal{B}_2^{\mathrm{rev}}.$$ We therefore find that $$\mathcal{Z}_{(2)}(\mathcal{A})=\mathcal{Z}_{(2)}(\mathcal{Z}(\mathcal{C},\mathcal{A}))\simeq \mathcal{Z}_{(2)}(\mathcal{B}_2^{\mathrm{rev}})\simeq \mathcal{E},$$ where the first equality arises from Lemma \ref{lem:symcentercentralizer} and Corollary \ref{cor:centernondeg}.
Consequently, it follows from Proposition \ref{prop:surjectivebraided} that the braided tensor functor $\mathcal{B}_1\rightarrow\mathcal{Z}(\mathcal{C})$ is necessarily fully faithful as desired, which concludes the proof of the theorem.
\end{proof}

\subsection{Absolutely Completely Anisotropic Braided Fusion 1-Categories}\label{sub:abscompaniso}

For use in the next section, we prove technical results about the Witt equivalence classes of braided fusion 1-categories over $\mathds{k}$. Recall that, in the theory of braided fusion 1-categories over an algebraically closed field $\overline{\mathds{k}}$, it was shown in \cite{DMNO} (see also \cite{DNO}) that every Witt equivalence class of a non-degenerate or slightly degenerate braided fusion 1-category has a unique representative given by a completely anisotropic braided fusion 1-category. For later use, we examine the generalization of this property to braided fusion 1-categories over an arbitrary field.

\begin{Definition}
A braided fusion 1-category $\mathcal{B}$ over $\mathds{k}$ is completely anisotropic if every commutative separable algebra $A$ in $\mathcal{B}$ is contained in the full fusion sub-1-category $\mathbf{Vect}_{\mathbb{K}}$ generated by the monoidal unit.
\end{Definition}

\noindent In particular, we emphasize that a completely anisotropic braided fusion 1-category is automatically either non-degenerate or slightly degenerate.

The above definition is not a priori well-adapted to the setting of non-algebraically closed fields. Instead, we will consider a strengthening of that notion. In what follows, we will implicitly use the fact established in Lemma \ref{lem:braidedcentrallinearity} that every braided fusion 1-category $\mathcal{B}$ over $\mathds{k}$ is automatically linear over $\mathbb{K}:=\Omega\mathcal{B}$.

\begin{Definition}
A braided fusion 1-category $\mathcal{B}$ is absolutely completely anisotropic if its base change $\overline{\mathcal{B}}$ from $\mathbb{K}$ to $\overline{\mathds{k}}$ is a completely anisotropic braided fusion 1-category.
\end{Definition}

\noindent Note that if $\mathcal{B}$ is absolutely completely anisotropic, then it must be completely anisotropic, and that, moreover, its base change $\mathcal{B}_{\mathbb{L}}$ for any (algebraic) field extension $\mathbb{L}/\mathbb{K}$ is again absolutely completely anisotropic.

We begin by proving that absolutely completely anisotropic braided fusion 1-categories exist.

\begin{Lemma}\label{lem:existenceabscompanis}
Let $\mathcal{B}$ be a non-degenerate or slightly degenerate braided fusion 1-category over $\mathds{k}$. There exists a finite field extension $\mathbb{L}/\mathbb{K}$ such that the base change $\mathcal{B}_{\mathbb{L}}$ is Witt equivalent to an absolutely completely anisotropic braided fusion 1-category.
\end{Lemma}
\begin{proof}
It was proven in \cite[Theorem 5.13]{DMNO} that $\overline{\mathcal{B}}$ is Witt equivalent to a completely anisotropic braided fusion 1-category over $\overline{\mathds{k}}$. More precisely, this means that there exists a $\overline{\mathcal{B}}$-central fusion 1-category $\overline{\mathcal{C}}$ over $\overline{\mathds{k}}$ for which the associated braided tensor functor $\overline{F}:\overline{\mathcal{B}}\hookrightarrow \mathcal{Z}(\overline{\mathcal{C}})$ is fully faithful and the braided fusion 1-category $\mathcal{Z}(\overline{\mathcal{C}},\overline{\mathcal{B}})$ over $\overline{\mathds{k}}$ is completely anisotropic. But, the fusion 1-category $\overline{\mathcal{C}}$ admits a form $\mathcal{C}_{\mathbb{L}}$ over a finite field extension $\mathbb{L}/\mathbb{K}$ (see for instance \cite[Remark 2.4]{ENO1}). Upon considering a larger finite field extension, thanks to Lemma \ref{lem:equivalencefinitebasechange}, we can also assume that the braided tensor functor $\overline{F}$ admits a form $F_{\mathbb{L}}:\mathcal{B}_{\mathbb{L}}\rightarrow \mathcal{Z}(\mathcal{C}_{\mathbb{L}})$ over $\mathbb{L}$. As $\overline{F}$ is fully faithful, so must be $F_{\mathbb{L}}$. Finally, the braided fusion 1-category $\mathcal{Z}(\mathcal{C}_{\mathbb{L}},\mathcal{B}_{\mathbb{L}})$ over $\mathbb{L}$ is absolutely completely anisotropic by construction and Lemma \ref{lem:centralizerbasechange}.
\end{proof}

Completely anisotropic braided fusion 1-categories, and a fortiori also the absolutely completely anisotropic ones, enjoy an incompressibility property that will play a key role below.

\begin{Lemma}\label{lem:centralincompressible}
Let $\mathcal{B}$ be a completely anisotropic braided fusion 1-category over $\mathds{k}$, and let $\mathcal{C}$ be a $\mathcal{B}$-central fusion 1-category. If $\Omega\mathcal{C} = \mathbb{K}$, then the canonical tensor functor $F:\mathcal{B}\rightarrow\mathcal{C}$ is fully faithful.
\end{Lemma}
\begin{proof}
Without loss of generality, we may assume that $\mathds{k}=\mathbb{K}$. Let $F^*$ be the right adjoint to $F$. Then, $A:=F^*(I)$ is a separable algebra in $\mathcal{B}$ by \cite[Corollary 2.6.9]{DSPS13}. The algebra $A$ is also commutative by the argument used in the proof of \cite[Lemma 3.5]{DMNO}. By construction, we also must have $End_{\mathcal{B}}(I,A)\cong \Omega\mathcal{C}=\mathbb{K}$. It follows from the fact that $\mathcal{B}$ is completely anisotropic that $A = I$, the monoidal unit of $\mathcal{B}$, so that $F$ must indeed be fully faithful by the argument used in the proof of \cite[Lemma 5.12]{DMNO}.
\end{proof}

We will use the above incompressibility property to prove a technical result that will be the key to the proof of our main theorem in the next section. It will make use of the following lemma concerning the Frobenius-Perron dimension as generalized in \cite[Definition 3.21]{Sa} to fusion 1-categories over an arbitrary field $\mathds{k}$. When $\mathds{k}=\overline{\mathds{k}}$ is algebraically closed, this recovers the classical notion of Frobenius-Perron dimension discussed for instance in \cite{EGNO}.

\begin{Lemma}\label{lem:FPdimGaloistrivial}
Let $\mathcal{C}$ be a fusion 1-category over $\mathds{k}$ with $\Omega\mathcal{C} = \mathds{k}$. Then, we have $\mathrm{FPdim}(\mathcal{C})=\mathrm{FPdim}(\overline{\mathcal{C}})$.
\end{Lemma}
\begin{proof}
For any simple object $X$ of $\mathcal C$, we have that $\mathbb D_X:=End_{\mathcal{C}}(X)$ is a division algebra of degree $n$ over its center $Z(\mathbb D_X)$, which is a field extension over $\mathds k$ of degree $k$.
Upon base extension to $\overline{\mathds k}$, $\mathbb D_X$ splits into $k$ copies of $M_n(\overline{\mathds k})$, which shows that $\overline{X}$, the image of $X$ under the canonical monoidal functor $\mathcal{C}\rightarrow \overline{\mathcal{C}}$, is isomorphic to $k$ many pairwise non-isomorphic simple objects $\overline{X}_i$, each with multiplicity $n$.
Since the absolute Galois group transitively permutes the $\overline{X}_i$, they must all have equal Frobenius-Perron dimension $\mathrm{FPdim}(\overline{X}_i)=d$.
Setting $D=\mathrm{FPdim}(X)$, we find that
\[k\cdot d^2\;\;=\;\;k\cdot\left(\frac{D}{nk}\right)^2\;\;=\;\;\frac{D^2}{\dim_{\mathbb K}(\mathbb D_X)}\;.\]
Using this formula, it follows from \cite[Thm 3.19]{Sa} that $\mathrm{FPdim}(\mathcal{C}) = \mathrm{FPdim}(\overline{\mathcal{C}})$.
\end{proof}

\begin{Proposition}\label{prop:technicalabsolutelycompletelyanisotropic}
Let $\mathcal{B}_1$ and $\mathcal{B}_2$ be two absolutely completely anisotropic braided fusion 1-category over $\mathds{k}$, and let $\mathcal{C}$ be a fusion 1-category witnessing a Witt equivalence over $\mathds{k}$ between them. Then, if $\Omega\mathcal{C} = \mathbb{K}$, we have $\mathcal{C}\simeq \mathcal{B}_1$ as $\mathcal{B}_1$-central fusion 1-categories and $\mathcal{B}_1\simeq \mathcal{B}_2$ as braided fusion 1-categories over $\mathds{k}$.
\end{Proposition}
\begin{proof}
From the equivalence $\mathcal{Z}(\mathcal{C},\mathcal{B}_1)\simeq \mathcal{B}_2^{\mathrm{rev}}$, it follows that $\mathcal{Z}(\mathcal{C},\mathcal{B}_1)$ is absolutely completely anisotropic.
We shall not use this last equivalence for the remainder of the proof, and may therefore safely assume that all of the functors that we will consider are linear over $\mathbb{K}$.
Namely, as $\mathcal{B}_1\hookrightarrow \mathcal{Z}(\mathcal{C})$ is fully faithful, we must have $\Omega\mathcal{Z}(\mathcal{C}) = \mathbb{K}$.
By hypothesis, $\Omega\mathcal{C} = \mathbb{K}$, so $\mathcal{C}$ must be Galois trivial.
Thanks to Lemma \ref{lem:FPdimGaloistrivial} above, we have that $\mathrm{FPdim}(\mathcal{C}) = \mathrm{FPdim}(\overline{\mathcal{C}})$, where $\overline{\mathcal{C}}$ denotes the base change of $\mathcal{C}$ from $\mathbb{K}$ to $\overline{\mathds{k}}$, and similarly for all the fusion 1-categories under consideration. In particular, by combining \cite[Theorem 3.14]{DGNO} with Lemma \ref{lem:centralizerbasechange}, we find that $$\mathrm{FPdim}(\mathcal{Z}(\mathcal{C})) = \mathrm{FPdim}(\mathcal{B}_1)\mathrm{FPdim}(\mathcal{Z}(\mathcal{C},\mathcal{B}_1))=\mathrm{FPdim}(\mathcal{B}_1)\mathrm{FPdim}(\mathcal{B}_2).$$ 
Now, by \cite[Theorem 4.9]{Sa}, we also have $\mathrm{FPdim}(\mathcal{Z}(\mathcal{C}))=\mathrm{FPdim}(\mathcal{C})^2$.
Moreover, by Lemma \ref{lem:centralincompressible}, the canonical tensor functors $\mathcal{B}_1\rightarrow \mathcal{C}$ and $\mathcal{Z}(\mathcal{C},\mathcal{B}_1)\rightarrow \mathcal{C}$ are fully faithful. It therefore follows that $\mathrm{FPdim}(\mathcal{C}) = \mathrm{FPdim}(\mathcal{B}_1)$, so that $\mathcal{B}_1\simeq \mathcal{C}$ as $\mathcal{B}_1$-central fusion 1-categories. Finally, we then have $$\mathcal{Z}(\mathcal{C},\mathcal{B}_1)\simeq\mathcal{Z}(\mathcal{B}_1,\mathcal{B}_1)\simeq \mathcal{B}_1^{\mathrm{rev}},$$ where the last equality comes from Corollary \ref{cor:cenralizernondegenerateincenter}.
This concludes the proof.
\end{proof}

Finally, for completeness, we prove a uniqueness statement for absolutely completely anisotropic braided fusion 1-categories that should be compared with \cite[Theorem 5.13]{DMNO}. 

\begin{Corollary}
Let $\mathcal{B}_1$ and $\mathcal{B}_2$ be two absolutely completely anisotropic braided fusion 1-categories over $\mathds{k}$. If $\mathcal{B}_1$ and $\mathcal{B}_2$ are Witt equivalent as braided fusion 1-categories over $\mathbb{K}$, then, for every sufficiently large finite field extension $\mathbb{L}/\mathbb{K}$, we have $(\mathcal{B}_1)_{\mathbb{L}}\simeq(\mathcal{B}_2)_{\mathbb{L}}$ as braided fusion 1-categories over $\mathbb{L}$.
\end{Corollary}
\begin{proof}
Let $\mathcal{C}$ be a fusion 1-category over $\mathbb{K}$ witnessing the Witt equivalence over $\mathbb{K}$ between $\mathcal{B}_1$ and $\mathcal{B}_2$. That is, we have a fully faithful braided tensor functor $\mathcal{B}_1\hookrightarrow \mathcal{Z}(\mathcal{C})$ and a braided equivalence $\mathcal{Z}(\mathcal{C},\mathcal{B}_1)\simeq \mathcal{B}_2^{\mathrm{rev}}$ over $\mathbb{K}$. Let $\mathbb{L}$ be a sufficiently large finite field extension of $\mathbb{K}$ such that $\Omega\mathcal{C}\otimes_{\mathbb{K}}\mathbb{L}$ splits as a direct sum of copies of $\mathbb{L}$. Note that, as $\Omega\mathcal{Z}(\mathcal{C}_{\mathbb{L}})\simeq \Omega \mathcal{Z}(\mathcal{C})_{\mathbb{L}}\simeq \mathbb{K}\otimes_{\mathbb{K}}\mathbb{L}$, the multifusion 1-category $\mathcal{C}_{\mathbb{L}}$ must be indecomposable. Upon replacing it with one of its fusion summands, we may therefore assume that $\mathcal{C}_{\mathbb{L}}$ is a fusion 1-category over $\mathbb{L}$.
As $\Omega\mathcal{C}_{\mathbb{L}}\simeq \mathbb{L}$ by construction, Proposition \ref{prop:technicalabsolutelycompletelyanisotropic} concludes the proof.
\end{proof}

\section{Compact Semisimple Tensor 2-Categories up to Morita Equivalence}\label{sec:CSS2CMorita}

Over an algebraically closed field of characteristic zero, it was shown in \cite{D9} that every fusion 2-category is Morita equivalent to a connected fusion 2-category. The proof can be summarized as follows:

\begin{enumerate}
    \item Every fusion 2-category is Morita equivalent to a fusion 2-category $\mathfrak{D}$ with $\Omega\mathfrak{D}$ non-degenerate or slightly degenerate.
    \item Every fusion 2-category $\mathfrak{D}$ such that $\Omega\mathfrak{D}$ is non-degenerate or slightly degenerate is Morita equivalent to a connected fusion 2-category.
\end{enumerate}

\noindent More precisely, the second step proceeds by showing that every fusion 2-category $\mathfrak{D}$ such that $\Omega\mathfrak{D}$ is non-degenerate or slightly degenerate is Morita equivalent to a fusion 2-category that has at least one invertible object in every connected component. The existence of these invertible objects is obtained by appealing to \cite{JFY}.

We will now work over an arbitrary field $\mathds{k}$ of characteristic zero.
As we explain in subsection \ref{sub:collapse} below, the first step in this argument can be readily carried out.
However, over an arbitrary field, a compact semisimple tensor 2-category $\mathfrak{D}$ such that $\Omega\mathfrak{D}$ is non-degenerate or slightly degenerate does not necessarily have an invertible object in every connected component.
Instead, we show in subsection \ref{sub:InflationHypernormal} that every such compact semisimple tensor 2-category is Morita equivalent to a compact semisimple tensor 2-category that has an invertible object in every connected component.
In order to do so, we employ a categorified version of the inflation procedure introduced in \cite{SS}.
We then carry out the analogue of the second step in section \ref{sub:Moritaconnected} using Theorem \ref{thm:Picardgrouptorsion}.

\subsection{First Step:\ Collapsing the Symmetric Center}\label{sub:collapse}

It is convenient to divide compact semisimple tensor 2-categories into two classes.

\begin{Definition}
A compact semisimple tensor 2-category is called bosonic if $\mathcal{Z}_{(2)}(\Omega\mathfrak{C})$ is Tannakian. Otherwise, it is called fermionic.
\end{Definition}

We have the following generalization of \cite[Theorem 4.1.6]{D9}.

\begin{Proposition}\label{prop:collapsingsymcenter}
Let $\mathfrak{C}$ be a bosonic, resp.\ fermionic, compact semisimple tensor 2-category over $\mathds{k}$.
Then, $\mathfrak{C}$ is Morita equivalent to a compact semisimple tensor 2-category $\mathfrak{D}$ with $\Omega\mathfrak{D}$ a non-degenerate, resp.\ slightly degenerate, braided fusion 1-category.
\end{Proposition}
\begin{proof}
For brevity, we only treat the bosonic case.
We write $\mathcal{B} = \Omega\mathfrak{C}$, and $\mathcal{E} = \mathcal{Z}_{(2)}(\mathcal{B})$.
By Theorem \ref{thm:Deligne}, $\mathcal{E}$ admits a fiber functor $F:\mathcal{E}\rightarrow\mathbf{Vect}_{\mathbb{L}}$, where $\mathbb{L}$ is a finite field extension of $\mathds{k}$.
The braided fusion 1-category $\mathcal{A} := \mathcal{B}\boxtimes_{\mathcal{E}}\mathbf{Vect}_{\mathbb{L}}$ has a canonical $\mathcal{B}$-central structure coming from the canonical surjective braided functor $\mathcal{B}\rightarrow \mathcal{A}$.
We may therefore view $\mathcal{A}$ as a connected rigid algebra in $\mathfrak{C}^0\subseteq \mathfrak{C}$.
The corresponding Morita dual compact semisimple tensor 2-category $\mathbf{Bimod}_{\mathfrak{C}}(\mathcal{A})$ satisfies $\Omega \mathbf{Bimod}_{\mathfrak{C}}(\mathcal{A}) = \mathcal{Z}(\mathcal{A},\mathcal{B})^{\mathrm{rev}}$ by Lemma \ref{lem:loopsdualconnected}.
Finally, observe that, by construction, we have that $\mathcal{Z}_{(2)}(\mathcal{A})\simeq \mathbf{Vect}_{\mathbb{L}}$, so that $\mathcal{Z}_{(2)}(\mathcal{Z}(\mathcal{A},\mathcal{B}))\simeq \mathbf{Vect}_{\mathbb{L}}$ by Lemma \ref{lem:symcentercentralizer}, concluding the proof.
\end{proof}

\begin{Example}
For any finite group $G$, we can consider the symmetric fusion 1-category $\mathbf{Rep}_{\mathds{k}}(G)$ of representations of $G$ in $\mathbf{Vect}_{\mathds{k}}$.
In particular, the forgetful functor is a symmetric monoidal functor $U:\mathbf{Rep}_{\mathds{k}}(G)\rightarrow \mathbf{Vect}_{\mathds{k}}$.
There is an equivalence $\mathbf{2Rep}_{\mathds{k}}(G)\simeq \mathbf{Mod}(\mathbf{Rep}_{\mathds{k}}(G))$ of compact semisimple tensor 2-categories.
Moreover, the forgetful 2-functor $\mathbf{2Rep}_{\mathds{k}}(G)\rightarrow \mathbf{2Vect}_{\mathds{k}}$ can be identified with $\mathbf{Mod}(U):\mathbf{Mod}(\mathbf{Rep}_{\mathds{k}}(G))\rightarrow \mathbf{Mod}(\mathbf{Vect}_{\mathds{k}})$.
It follows that the Morita dual to $\mathbf{Mod}(\mathbf{Rep}_{\mathds{k}}(G))$ with respect to $\mathbf{Mod}(\mathbf{Vect}_{\mathds{k}})$ is $\mathbf{2Vect}_{\mathds{k}}(G)$, the compact semisimple tensor 2-category of $G$-graded 2-vector spaces over $\mathds{k}$.
In particular, it has an invertible object in every connected component.
\end{Example}

\begin{Example}
For any finite abelian group $A$, we can consider the symmetric fusion 1-category $\mathbf{Vect}_{\mathds{k}}(A)$ of $A$-graded (finite) vector spaces over $\mathds{k}$.
The forgetful functor $U:\mathbf{Vect}_{\mathds{k}}(A)\rightarrow \mathbf{Vect}_{\mathds{k}}$ is symmetric monoidal.
There is an equivalence $\mathbf{Mod}(\mathbf{Vect}_{\mathds{k}}(A))\simeq \mathbf{2Vect}_{\mathds{k}}(A[1])$ of compact semisimple tensor 2-categories.
Here, $A[1]$ denotes the 2-group with a single object and $A$ as its group of automorphisms, and the compact semisimple tensor 2-category $\mathbf{2Vect}_{\mathds{k}}(A[1])$ is defined as in \cite[Section 2.1.3]{DR}.
As in the previous example, the forgetful functor $\mathbf{2Vect}_{\mathds{k}}(A[1])\rightarrow \mathbf{2Vect}_{\mathds{k}}$ can be identified with $\mathbf{Mod}(U)$.
The Morita dual compact semisimple tensor 2-category to $\mathbf{Mod}(\mathbf{Vect}_{\mathds{k}}(A))$ with respect to $\mathbf{Mod}(\mathbf{Vect}_{\mathds{k}})$ is $\mathbf{2Rep}_{\mathds{k}}(A[1])$.
If $\mathds{k}$ has enough roots of unity, so that $Hom(A,\mathds{k}^\times)\cong A$ non-canonically, then $\mathbf{2Rep}_{\mathds{k}}(A[1])\simeq \mathbf{2Vect}_{\mathds{k}}(A)$ by the argument of \cite[Example 1.4.22]{DR}.
However, in general, this is not the case, and $\mathbf{2Rep}_{\mathds{k}}(A[1])$ does not contain an invertible object in every connected component (see for instance \cite[Remark 3.6]{D11}, which considers the case $G=\mathbb{Z}/3\mathbb{Z}$ and $\mathds{k}=\mathbb{R}$).
\end{Example}

We will need the following slight refinement of Proposition \ref{prop:collapsingsymcenter} above.

\begin{Corollary}
Let $\mathfrak{C}$ be a compact semisimple tensor 2-category over $\mathds{k}$.
Then, $\mathfrak{C}$ is Morita equivalent to a compact semisimple tensor 2-category $\mathfrak{D}$ with $\Omega\mathfrak{D}$ an absolutely completely anisotropic braided fusion 1-category.
\end{Corollary}
\begin{proof}
By Proposition \ref{prop:collapsingsymcenter}, it is enough to consider the case when $\mathcal{B}:=\Omega\mathfrak{C}$ is non-degenerate or slightly degenerate. Thanks to Lemma \ref{lem:existenceabscompanis}, there exists a finite field extension $\mathbb{L}/\mathbb{K}$ such that the base change $\mathcal{B}_{\mathbb{L}}$ is Witt equivalent to an absolutely completely anisotropic braided fusion 1-category, that is, there exists a $\mathcal{B}_{\mathbb{L}}$-central fusion 1-category $\mathcal{C}$ such that $\mathcal{Z}(\mathcal{C},\mathcal{B}_{\mathbb{L}})$ is absolutely completely anisotorpic. 
Let us note that we also have $\mathcal{Z}(\mathcal{C},\mathcal{B})=\mathcal{Z}(\mathcal{C},\mathcal{B}_{\mathbb{L}})$ as the dominant (full) image of $\mathcal{B}$ in $\mathcal{Z}(\mathcal{C})$ is $\mathcal{B}_{\mathbb{L}}$.
We may view $\mathcal{C}$ as a $\mathcal{B}$-central fusion 1-category, and therefore also as a connected rigid algebra in $\mathfrak{C}^0=\mathbf{Mod}(\mathcal{B})$, so that we can set $\mathfrak{D}:=\mathbf{Bimod}_{\mathfrak{C}}(\mathcal{C})$. Thanks to Lemma \ref{lem:loopsdualconnected}, we find $$\Omega\mathfrak{D} \simeq \mathcal{Z}(\mathcal{C},\mathcal{B})^{\mathrm{rev}}=\mathcal{Z}(\mathcal{C},\mathcal{B}_{\mathbb{L}})^{\mathrm{rev}},$$ which is absolutely completely anisotropic by construction.
\end{proof}

\subsection{Second Step:\ Inflation}\label{sub:InflationHypernormal}

Our goal is to establish the following result.

\begin{Theorem}\label{thm:MoritaGroupTheoretical}
Let $\mathfrak{C}$ be a compact semisimple tensor 2-category with $\mathcal{B}:=\Omega\mathfrak{C}$ an absolutely completely anisotropic braided fusion 1-category.
Then, $\mathfrak{C}$ is Morita equivalent to a compact semisimple tensor 2-category that has an invertible object in every connected component.
\end{Theorem}

Inspired by the study of Morita equivalence classes of fusion 1-category over $\mathds{k}$ undergone in \cite{SS}, it is natural to suspect that the procedure of inflation will play a key role.

\begin{Construction}\label{con:inflation}
Let $\mathfrak{C}$ be a compact semisimple tensor 2-category, and write $\mathcal{B}:=\Omega\mathfrak{C}$.
Let $\mathbb{L}$ be a (finite normal) field extension of $\mathbb{K}:=\Omega\mathcal{B}$.
We may then consider the braided fusion 1-category $\mathcal{B}_{\mathbb{L}}$ as a $\mathcal{B}$-central fusion 1-category.
In particular, in the language of \cite{D7}, $\mathcal{B}_{\mathbb{L}}$ can be viewed as a connected rigid algebra in $\mathfrak{C}^0=\mathbf{Mod}(\mathcal{B})$.
Note that, by definition, $\Omega\mathcal{B}_{\mathbb{L}} = \mathbb{L}$.
We define the inflation of $\mathfrak{C}$ from $\mathbb{K}$ to $\mathbb{L}$ by $$\mathrm{Inf}_{\mathbb{K}}^{\mathbb{L}}(\mathfrak{C}):=\mathbf{Bimod}_{\mathfrak{C}}(\mathcal{B}_{\mathbb{L}}).$$ By construction, the inflation is a compact semisimple tensor 2-category that is Morita equivalent to $\mathfrak{C}$.
Below, we record some elementary properties of the inflation $\mathrm{Inf}_{\mathbb{K}}^{\mathbb{L}}(\mathfrak{C})$.
\end{Construction}

\begin{Lemma}\label{lem:inflationbasic}
Let $\mathbb{L}$ be a (finite normal) field extension of $\mathbb{K}$.
Then, the compact semisimple tensor 2-category $\mathrm{Inf}_{\mathbb{K}}^{\mathbb{L}}(\mathfrak{C})$ satisfies $\Omega \mathrm{Inf}_{\mathbb{K}}^{\mathbb{L}}(\mathfrak{C})\simeq \mathcal{B}_{\mathbb{L}}$.
In particular, if $\mathcal{B}=\Omega\mathfrak{C}$ is absolutely completely anisotropic then so is $\Omega\mathrm{Inf}_{\mathbb{K}}^{\mathbb{L}}(\mathfrak{C})$.
\end{Lemma}
\begin{proof}
The full image of $\mathcal{B}$ in $\mathcal{Z}(\mathcal{B}_{\mathbb{L}})$ can be identified with $\mathcal{B}_{\mathbb{L}}$.
By Lemma \ref{lem:loopsdualconnected}, it is therefore enough to show that $\mathcal{Z}(\mathcal{B}_{\mathbb{L}},\mathcal{B}_{\mathbb{L}})\simeq \mathcal{B}_{\mathbb{L}}^{\mathrm{rev}}$, which was done in Corollary \ref{cor:cenralizernondegenerateincenter}.
\end{proof}

We now introduce a strengthening of the notion of normality introduced in definition \ref{def:normal} above.

\begin{Definition}\label{def:hypernomraltensor}
A simple object $C$ of $\mathfrak{C}$ is hypernormal if $\Omega End_{\mathfrak{C}}(C)\cong\Omega^2\mathfrak{C}$, and a component $[C]$ of $\mathfrak{C}$ is said to be hypernormal if it contains a hypernormal simple object.
The compact semisimple tensor 2-category $\mathfrak{C}$ is hypernormal if all of its connected components are hypernormal.
\end{Definition}

\noindent Observe that every hypernormal simple object $C$ is necessarily normal.
Namely, $\Omega \mathcal{Z}(End_{\mathfrak{C}}(C))$ must contain $\mathbb{K}=\Omega^2\mathfrak{C}$ as $\mathfrak{C}^0$ acts on the connected component of $C$, but it also cannot be any larger by hypernormality of $C$.
We emphasize that, unlike normality, it is not true that every simple object of a hypernormal connected component is hypernormal.

We use inflation to produce a Morita equivalent compact semisimple tensor 2-category that is hypernormal.

\begin{Lemma}\label{lem:inflationhypernormal}
Let $\mathfrak{C}$ be a compact semisimple tensor 2-category.
For every large enough finite normal field extension $\mathbb{L}$ of $\mathbb{K}=\Omega^2\mathfrak{C}$, the inflation $\mathrm{Inf}_{\mathbb{K}}^{\mathbb{L}}(\mathfrak{C})$ is hypernormal.
\end{Lemma}
\begin{proof}
Let $\mathcal{C}$ be a multifusion 1-category over $\mathds{k}$ such that $\mathfrak{C}\simeq \mathbf{Mod}(\mathcal{C})$ as compact semisimple 2-categories.
In fact, as $\mathfrak{C}\simeq \mathbf{Mod}(\mathcal{C})$ is a module 2-category over the compact semisimple tensor 2-category $\mathfrak{C}^0\simeq \mathbf{Mod}(\mathcal{B})$, it follows that $\mathcal{C}$ is canonically a $\mathcal{B}$-central multifusion 1-category.
In particular, this endows $\mathcal{C}$ with a $\mathbb{K}$-linear structure.
We say that the finite normal field extension $\mathbb{L}/\mathbb{K}$ is large enough if $\Omega\mathcal{C}_{\mathbb{L}} = \Omega\mathcal{C}\otimes_{\mathbb{K}}\mathbb{L}$ splits as a direct sums of copies of $\mathbb{L}$.

By definition, we have $\mathrm{Inf}_{\mathbb{K}}^{\mathbb{L}}(\mathfrak{C})=\mathbf{Bimod}_{\mathfrak{C}}(\mathcal{B}_{\mathbb{L}})$, so that it is enough to show that every connected component of $\mathbf{Bimod}_{\mathfrak{C}}(\mathcal{B}_{\mathbb{L}})$ is hypernormal.
We have equivalences of compact semisimple 2-categories: \begin{equation}\label{eq:inflationunderlying2cat}\mathbf{Bimod}_{\mathfrak{C}}(\mathcal{B}_{\mathbb{L}})\simeq \mathbf{Mod}(\mathcal{B}^{\mathrm{mop}}_{\mathbb{L}}\boxtimes_{\mathcal{B}}\mathcal{C}\boxtimes_{\mathcal{B}}\mathcal{B}_{\mathbb{L}})\simeq \mathbf{Mod}(\mathbf{Vect}_{\mathbb{L}}\boxtimes_{\mathbb{K}}\mathcal{C}\boxtimes_{\mathbb{K}}\mathbf{Vect}_{\mathbb{L}}),\end{equation} using that $\mathcal{B}_{\mathbb{L}}=\mathcal{B}\boxtimes_{\mathbb{K}}\mathbf{Vect}_{\mathbb{L}}$.
But, we have that $\Omega(\mathbf{Vect}_{\mathbb{L}}\boxtimes_{\mathbb{K}}\mathcal{C}\boxtimes_{\mathbb{K}}\mathbf{Vect}_{\mathbb{L}})=\mathbb{L}\otimes_{\mathbb{K}}\Omega\mathcal{C}\otimes_{\mathbb{K}}\mathbb{L}$ by inspection, which is a direct sums of copies of $\mathbb{L}$ by our choice of $\mathbb{L}$.
This shows that $\mathbf{Bimod}_{\mathfrak{C}}(\mathcal{B}_{\mathbb{L}})$ is hypernormal and thereby concludes the proof of the lemma.
\end{proof}

We are now ready to move towards the proof of Theorem \ref{thm:MoritaGroupTheoretical}.
Thanks to the last lemma above, it is enough to consider compact semisimple tensor 2-categories $\mathfrak{C}$ for which not only the braided fusion 1-categories $\Omega\mathfrak{C}$ are absolutely completely anisotropic, and therefore either non-degenerate or slightly degenerate, but are also hypernormal.
In particular, such a compact semisimple tensor 2-category is normal and therefore admits a Galois grading by Proposition \ref{prop:Galoisgrading}.
We begin by considering the case when the aforementioned grading is trivial.

\begin{Lemma}\label{lem:normaltrivialGalois}
Let $\mathfrak{C}$ be a normal compact semisimple tensor 2-category over $\mathds{k}$ such that $\mathcal{B}=\Omega\mathfrak{C}$ is either non-degenerate or slightly degenerate.
If the Galois grading on $\mathfrak{C}$ is trivial, then the connected components of $\mathfrak{C}$ inherit a group structure from the monoidal product, and therefore supply $\mathfrak{C}$ with a faithful grading.
\end{Lemma}
\begin{proof}
As the Galois grading on $\mathfrak{C}$ is trivial by assumption, the compact semisimple tensor 2-category $\mathfrak{C}$ is linear over $\Omega^2\mathfrak{C}$.
We may therefore without loss of generality assume that $\mathds{k}=\Omega^2\mathfrak{C}$.
Let us write $\overline{\mathfrak{C}}$ for the base change of $\mathfrak{C}$ to the algebraic closure $\overline{\mathds{k}}$ of $\mathds{k}$.
In particular, there is a monoidal 2-functor $F:\mathfrak{C}\rightarrow \overline{\mathfrak{C}}$.
Since $\Omega^2\mathfrak{C}=\mathds{k}$, it follows that $\Omega^2\overline{\mathfrak{C}}=\overline{\mathds{k}}$, so that $\overline{\mathfrak{C}}$ is a fusion 2-category over $\overline{\mathds{k}}$.
In fact, we can say more:\ As every connected component of $\mathfrak{C}$ is normal, it follows that the monoidal 2-functor $F$ must induce a bijection $\pi_0(\mathfrak{C})\cong \pi_0(\overline{\mathfrak{C}})$ on connected components.
Namely, if $C$ is a simple object of $\mathfrak{C}$, then $\Omega\mathcal{Z}(End_{\mathfrak{C}}(C)) = \mathds{k}$ by normality.
But then, we have $\Omega\mathcal{Z}(End_{\overline{\mathfrak{C}}}(F(C))) = \overline{\mathds{k}}$, so that $F(C)$ is contained in a single connected component of $\overline{\mathfrak{C}}$.
Finally, it follows from Lemma \ref{lem:centralizerbasechange} that $\Omega\overline{\mathfrak{C}}$ is non-degenerate or slightly degenerate.
In either case, it follows from Lemma \ref{lem:connectedcomponentsgroup} that the finite set $\pi_0(\overline{\mathfrak{C}})$ inherits a group structure $G$ from the monoidal product of $\overline{\mathfrak{C}}$. More precisely, the fusion 2-category $\overline{\mathfrak{C}}$ is a quasi-trivial faithfully $G$-graded extension of $\overline{\mathfrak{C}}^0$.
Putting the above discussion together, we find that $\mathfrak{C}$ is a faithfully $G$-graded extension of $\mathfrak{C}^0$.
\end{proof}

\begin{Proposition}\label{prop:hypernormalgradedcomponent}
Let $\mathfrak{C}$ be a hypernormal compact semisimple tensor 2-category over $\mathds{k}$ such that $\mathcal{B}=\Omega\mathfrak{C}$ is absolutely completely anisotropic.
Then, the connected components of $\mathfrak{C}$ inherit a group structure from the monoidal product, and therefore supply $\mathfrak{C}$ with a faithful grading.
\end{Proposition}
\begin{proof}
As $\mathfrak{C}$ is normal, it sports a canonical Galois grading by $\mathrm{Gal}(\mathbb{K}/\mathds{k})$.
We will write $\mathfrak{C}_e$ for the trivially Galois graded part.
Let us also write $\Gamma\subseteq \mathrm{Gal}(\mathbb{K}/\mathds{k})$ for the subgroup corresponding to the support of $\mathfrak{C}$, i.e.\ we have $\mathfrak{C} = \boxplus_{\gamma\in\Gamma}\mathfrak{C}_{\gamma}$ with $\mathfrak{C}_{\gamma}$ non-zero.
We begin by observing that it will be enough to show that, for every simple object $C$ of $\mathfrak{C}$, the monoidal product $C^{\sharp}\Box C$ of $C$ with its dual lies entirely in the connected component of the identity $\mathfrak{C}^0$.
Namely, provided that this is the case, then, for every simple objects $C$ and $D$ in $\mathfrak{C}$, we have that $C^{\sharp}\Box C\Box D$ lies in a single connected component.
This is possible only if $C\Box D$ lies in a single connected component, so that $\mathfrak{C}$ is indeed faithfully graded by its connected components.

It remains to prove that $\mathfrak{C}$ satisfies the above property. Let us fix an arbitrary element $\gamma$ of $\Gamma$.
Thanks to Lemma \ref{lem:normaltrivialGalois}, the compact semisimple tensor 2-category $\mathfrak{C}_e$ is graded by its connected components.
Moreover, it follows from \cite[Theorem 5.3.4]{D4} that there exists a connected rigid algebra $A$ in $\mathfrak{C}_e$ such that $\mathfrak{C}_{\gamma}\simeq\mathbf{Mod}_{\mathfrak{C}_e}(A)$ as left $\mathfrak{C}_e$-module 2-categories.
Let us record that by hypernormality it follows from \cite[Section 5.3]{D4} that we can choose $A$ so that $End_{A}(Id_A)\cong \mathbb{K}$.
As $\mathfrak{C}_{\gamma}\simeq\mathbf{Mod}_{\mathfrak{C}_e}(A)$ is an invertible $\mathfrak{C}_e$-bimodule 2-category, we have that $\Omega\mathbf{Bimod}_{\mathfrak{C}_e}(A)\simeq\Omega\mathfrak C_e=\mathcal{B}$ as braided fusion 1-categories over $\mathds{k}$ by definition, so that $\Omega\mathbf{Bimod}_{\mathfrak{C}_e}(A)$ is an absolutely completely anisotropic braided fusion 1-category over $\mathbb{K}$. 
The normal compact semisimple tensor 2-category $\mathfrak{C}_e$ is linear over $\mathbb{K}$ by triviality of its Galois grading, we can therefore consider its base change $\overline{\mathfrak{C}_e}$ from $\mathbb{K}$ to its algebraic closure. 
There is a monoidal 2-functor $\mathfrak{C}_e\rightarrow \overline{\mathfrak{C}_e}$, and, abusing notations slightly, we also write $A$ for the image of the algebra $A$ in $\mathfrak{C}_e$ under this monoidal 2-functor.
The algebra $A$ in $\overline{\mathfrak{C}_e}$ is not only rigid, since this property is preserved by every monoidal 2-functor, but also connected because the algebra $A$ in $\mathfrak{C}_e$ satisfies $End_{A}(Id_A)\cong \mathbb{K}$.
It was established in \cite[Proof of Thm.\ 3.1.4]{D10} that the base changes of $\mathbf{Bimod}_{\mathfrak{C}_e}(A)$ and $\mathbf{Mod}_{\mathfrak{C}_e}(A)$ from $\mathbb{K}$ to its algebraic closure can be described as $$\overline{\mathbf{Bimod}_{\mathfrak{C}_e}(A)}\simeq \mathbf{Bimod}_{\overline{\mathfrak{C}_e}}(A)\,,\ \overline{\mathbf{Mod}_{\mathfrak{C}_e}(A)}\simeq \mathbf{Mod}_{\overline{\mathfrak{C}_e}}(A).$$
Now, because $\mathbf{Bimod}_{\mathfrak{C}_e}(A)$ and $\mathbf{Mod}_{\mathfrak{C}_e}(A)$ are hypernormal, as they can be identified with full sub-2-categories of $\mathfrak{C}$, we have that the canonical 2-functors $\mathbf{Bimod}_{\mathfrak{C}_e}(A)\rightarrow \mathbf{Bimod}_{\overline{\mathfrak{C}_e}}(A)$ and $\mathbf{Mod}_{\mathfrak{C}_e}(A)\rightarrow \mathbf{Mod}_{\overline{\mathfrak{C}_e}}(A)$ induce bijections on connected components.
In particular, it follows from Lemma \ref{lem:Moritageneralizedstronglyfusion} that not only must $A$ be contained in $\mathfrak{C}^0_e$, but also that the forgetful 2-functors $\mathbf{Bimod}_{\mathfrak{C}_e}(A)\rightarrow \mathfrak{C}_e$ and $\mathbf{Mod}_{\mathfrak{C}_e}(A)\rightarrow\mathfrak{C}_e$ induce bijections on connected components. We will make use of this property below.

Finally, let $C$ be a simple object of $\mathfrak{C}_{\gamma}$. Its dual $C^{\sharp}$ is a simple object in $\mathfrak{C}_{\gamma^{-1}}$. It was shown in \cite[Proposition 3.4]{D11} that the 2-functor $$\Box:\mathfrak{C}_{\gamma^{-1}}\times\mathfrak{C}_{\gamma}\rightarrow \mathfrak{C}_e$$ can be identified with the canonical 2-functor $$\Box: \mathbf{LMod}_{\mathfrak{C}_e}(A)\times \mathbf{Mod}_{\mathfrak{C}_e}(A)\rightarrow \mathbf{Bimod}_{\mathfrak{C}_e}(A).$$ Moreover, there is a commutative diagram $$\begin{tikzcd}[sep=small]
\mathbf{LMod}_{\mathfrak{C}_e}(A)\times\mathbf{Mod}_{\mathfrak{C}_e}(A) \arrow[rr, "\Box"] \arrow[d] &  & \mathbf{Bimod}_{\mathfrak{C}_e}(A) \arrow[d] \\
\mathfrak{C}_e\times\mathfrak{C}_e \arrow[rr, "\Box"']&&\mathfrak{C}_e,
\end{tikzcd}$$ in which the vertical arrows are forgetful 2-functors. Now, it was established in the preceding paragraph that the vertical arrows induce bijections on the sets of connected components.
Furthermore, as $\mathfrak{C}_e$ is graded by its connected components thanks to Lemma \ref{lem:normaltrivialGalois}, the monoidal product of any two simple objects of $\mathfrak{C}_e$ must lie in a single connected component. This shows that the tensor product $C^{\sharp}\Box C$ must lie in a single connected component of $\mathfrak{C}_e$ as claimed.
As $\gamma$ and $C$ were arbitrary, this concludes the proof.
\end{proof}

Combining the above results, we find that, for every sufficiently large finite normal field extension $\mathbb{L}$ of $\mathbb{K}$, the compact semisimple tensor 2-category $\mathrm{Inf}_{\mathbb{K}}^{\mathbb{L}}(\mathfrak{C})$ is graded by its connected components.
It will follow from the fact that we have assumed that $\mathcal{B}=\Omega\mathfrak{C}$ is absolutely completely anisotropic that $\mathrm{Inf}_{\mathbb{K}}^{\mathbb{L}}(\mathfrak{C})$ has an invertible object in every connected component. 

\begin{Lemma}\label{lem:invertibleconnectednondeg}
Let $\mathcal{B}$ be an absolutely completely anisotropic braided fusion 1-category over $\mathds{k}$ with $\mathbb{K}:=\Omega\mathcal{B}$, and let $\mathfrak{M}$ be an invertible $\mathbf{Mod}(\mathcal{B})$-bimodule 2-category.
If $\Omega End_{\mathfrak{M}}(M) = \mathbb{K}$ for any simple object $M$ of $\mathfrak{M}$, then $\mathfrak{M}$ is a quasi-trivial $\mathbf{Mod}(\mathcal{B})$-bimodule 2-category.
\end{Lemma}
\begin{proof}
It follows from Theorem \ref{thm:MoritaWitt} that equivalence classes of invertible $\mathbf{Mod}(\mathcal{B})$-bimodule 2-categories are in bijection with Morita equivalence classes of fusion 1-categories $\mathcal{C}$ equipped with both a fully faithful braided tensor functor $\mathcal{B}\hookrightarrow \mathcal{Z}(\mathcal{C})$ and an equivalence $\mathcal{Z}(\mathcal{C},\mathcal{B})\simeq \mathcal{B}^{\mathrm{rev}}$ of braided fusion 1-categories over $\mathds{k}$.
The result is therefore a reformulation of Proposition \ref{prop:technicalabsolutelycompletelyanisotropic} above.
\end{proof}

We are finally in a position to prove Theorem \ref{thm:MoritaGroupTheoretical}.

\begin{proof}[Proof of Thm.\ref{thm:MoritaGroupTheoretical}]
It follows from Lemmas \ref{lem:inflationbasic} and \ref{lem:inflationhypernormal} that we may assume that $\mathfrak{C}$ is hypernormal.
In particular, it then follows from Proposition \ref{prop:hypernormalgradedcomponent} that the monoidal product induces a group structure $G$ on $\pi_0(\mathfrak{C})$, the set of connected component of $\mathfrak{C}$, so that $\mathfrak{C}$ is a faithfully $G$-graded extension of $\mathfrak{C}^0$.
By assumption, we have that $\mathfrak{C}^0=\mathbf{Mod}(\mathcal{B})$ with $\mathcal{B}$ absolutely completely anisotropic. But, thanks to Lemma \ref{lem:invertibleconnectednondeg}, any group graded extension of $\mathfrak{C}^0=\mathbf{Mod}(\mathcal{B})$ is quasi-trivial. Lemma \ref{lem:quasitrivialinvertible} therefore concludes the proof of the theorem.
\end{proof}

\subsection{Third Step:\ Morita Connectedness}\label{sub:Moritaconnected}

Provided that a compact semisimple tensor 2-category contains (at least) one invertible object in every connected component, then one can hope to adapt the results of \cite[Section 4.2]{D9} to show that it is Morita equivalent to a connected one.
While this may seem straightforward at first glance, a technical complication arises when working over an arbitrary field $\mathds{k}$ (of characteristic zero):\ A compact semisimple tensor 2-category can have infinitely many invertible objects.
However, this difficulty can be addressed using Theorem \ref{thm:Picardgrouptorsion}.

\begin{Theorem}\label{thm:MoritaConnected}
Let $\mathfrak{C}$ be compact semisimple 2-category that has (at least) one invertible object in every connected component, then it is Morita equivalent to a connected compact semisimple tensor 2-category.
\end{Theorem}
\begin{proof}
Let $\mathcal{B}:=\Omega\mathfrak{C}$, and $\mathbb{K}:=\Omega\mathcal{B}=\Omega^2\mathfrak{C}$.
Recall that an algebra is called strongly connected if its unit 1-morphism $u:I\rightarrow A$ is the inclusion of a simple summand.
Thanks to Lemma \ref{lem:EnoughInvertibleMoritaConnected} above, it is enough to exhibit a strongly connected rigid algebra $A$ in $\mathfrak{C}$ whose underlying object has a summand in every connected component of $\mathfrak{C}$.

To this end, note that the group of invertible objects of $\mathfrak{C}$ fits into a short exact sequence of groups $$0\rightarrow \mathrm{Pic}(\mathcal{B})\rightarrow \mathrm{Inv}(\mathfrak{C})\rightarrow \pi_0(\mathfrak{C})\rightarrow 1,$$ where the group $\pi_0(\mathfrak{C})$ is finite and the group $\mathrm{Pic}(\mathcal{B})$ is indfinite thanks to Theorem \ref{thm:Picardgrouptorsion}. It therefore follows from Lemma \ref{lem:indfiniteextension} that $\mathrm{Inv}(\mathfrak{C})$ is indfinite, so that there exists a finite subgroup $\widehat{G}$ of $\mathrm{Inv}(\mathfrak{C})$ that surjects onto $\pi_0(\mathfrak{C})$.
Then, we can consider the monoidal sub-2-category $\mathfrak{C}^{\times}$ of $\mathfrak{C}$ on the invertible objects and invertible morphisms, as well as its monoidal full sub-2-category $\mathfrak{G} \subseteq \mathfrak{C}^{\times}$ on the objects in $\widehat{G}$.
By construction, the first two homotopy groups of $\mathfrak{G}$ are finite, and its top homotopy group is $\mathbb{K}^{\times}$, which potentially sports a non-trivial action by $\widehat{G}$. The group-like topological monoid $\mathfrak{G}$ admits a Postnikov decomposition $\mathcal{G}\times^{\pi}\mathrm{B}^2\mathbb{K}^{\times}$ as the extension of the finite 2-group $\mathcal{G}$ by $\mathrm{B}^2\mathbb{K}^{\times}$ classified by a 4-cocycle $\pi$ in $H^4(\mathrm{B}\mathcal{G};\mathbb{K}^{\times})$.
Now, by Lemma \ref{lem:technicalcocycletrivialization} below, there exists a finite group $\widetilde{G}$ and an essentially surjective monoidal functor $p:\widetilde{G}\rightarrow \mathcal{G}$ such that the pullback $p^*\pi$ is trivializable.
In particular, we have a monoidal 2-functor $\widetilde{G}\rightarrow \mathfrak{C}$ whose image contains one invertible object in every connected component of $\mathfrak{C}$. We therefore also have a monoidal 2-functor $\mathrm{B}^2\mathbb{K}\rtimes\widetilde{G}\rightarrow \mathfrak{C}$, where $\widetilde{G}$ acts on $\mathbb{K}$ via the group homomorphism $\widetilde{G}\rightarrow \widehat{G}$ induced by $p$.
Upon taking Cauchy completion, we obtain a monoidal 2-functor $\mathbf{2Vect}_{\mathbb{K}}(\widetilde{G})\rightarrow \mathfrak{C}$ that hits every connected component of the target.
The image of the strongly connected rigid algebra $\mathbf{Vect}_{\mathbb{K}}(\widetilde{G})$ in $\mathbf{2Vect}_{\mathbb{K}}(\widetilde{G})$ is then a strongly connected rigid algebra in $\mathfrak{C}$ whose underlying object has a summand in every connected component of $\mathfrak{C}$ by construction.
This concludes the proof.
\end{proof}

The following technical result, which is a slight generalization of \cite[Lemma 4.2.4]{D9} was used in the proof of the above theorem.

\begin{Lemma}\label{lem:technicalcocycletrivialization}
Let $\mathcal{G}$ be any finite 2-group, $A$ an arbitrary abelian group (potentially equipped with an action of $\mathcal{G}$), and $\pi$ any 4-cocycle for $\mathcal{G}$ with coefficients in $A$. Then, there exists a finite group $\widetilde{G}$ and an essentially surjective monoidal functor $p:\widetilde{G}\rightarrow \mathcal{G}$ such that the pullback of $\pi$ along $p$ is a coboundary.
\end{Lemma}
\begin{proof}
The proof is almost exactly that given in \cite[Lemma 4.2.4]{D9} except that we need to appeal to a generalization of \cite[Lemma 4.2.3]{D9}. More precisely, we need to use the following statement:\ Let $G$ be a finite group, and $A$ an arbitrary abelian group equipped with an action by $G$.
For $n\geq 3$, given any $n$-cocycle $\omega$ for $G$ with coefficients in $A$, there exists a finite group $\widetilde{G}$ together with a surjective group homomorphism $p:\widetilde{G}\rightarrow G$ such that the pullback $p^*\omega$ is a coboundary.
Under the additional assumption that $A$ is torsion, this is exactly the statement of \cite[Lemma 4.2.3]{D9}.
However, upon inspection, we find that the proof only uses that the group $H^{n-2}(G;A)$ is torsion.
As we have assumed that $n\geq 3$, the group $H^{n-2}(G;A)$ is torsion for any arbitrary abelian group $A$, which concludes the proof.
\end{proof}

\section{Applications}\label{sec:applications}

We give some applications of our main theorems to the study of Witt groups of braided fusion 1-categories over $\mathds{k}$.

\subsection{Witt Groups}

The Witt group $\mathrm{Witt}(\mathbf{Vect}_{\mathds{k}})$ is the group whose elements are Witt equivalence classes of braided fusion 1-categories whose symmetric center is $\mathbf{Vect}_{\mathds{k}}$. This group admits a natural higher categorical interpretation. A symmetric monoidal 4-category $\mathbf{BrFus}_{\mathds{k}}$ whose objects are braided multifusion 1-categories over $\mathds{k}$ and whose 1-morphisms are central multifusion 1-categories was constructed in \cite{BJS}. Moreover, its invertible objects are exactly the braided fusion 1-categories whose symmetric center is $\mathbf{Vect}_{\mathds{k}}$ by \cite{BJSS}, and the group $\mathrm{Witt}(\mathbf{Vect}_{\mathds{k}})$ coincides with the group of equivalence classes of invertible objects of $\mathbf{BrFus}_{\mathds{k}}$. By \cite{D10}, there is also a symmetric monoidal 4-category $\mathbf{CSST2C}_{\mathds{k}}$ of compact semisimple tensor 2-categories over $\mathds{k}$, whose monoidal structure is given by the 2-Deligne tensor product of \cite{D3}. Its invertible objects can be characterized as follows by \cite[Lemma 4.3.5]{D11}.

\begin{Lemma}
A compact semisimple tensor 2-category $\mathfrak{C}$ over $\mathds{k}$ is (Morita) invertible, i.e.\ $\mathfrak{C}\boxtimes_{\mathds{k}}\mathfrak{C}^{\mathrm{mop}}$ is Morita equivalent to $\mathbf{2Vect}_{\mathds{k}}$, if and only if $\mathscr{Z}(\mathfrak{C})\simeq \mathbf{2Vect}_{\mathds{k}}$ as compact semisimple tensor 2-categories.
\end{Lemma}

\noindent Combining Theorems \ref{maintheorem} and \ref{thm:MoritaWitt}, we immediately obtain the following statement.

\begin{Corollary}
The Witt group $\mathrm{Witt}(\mathbf{Vect}_{\mathds{k}})$ of braided fusion 1-categories whose symmetric center is $\mathbf{Vect}_{\mathds{k}}$ is isomorphic to the group of (Morita equivalence classes of) invertible compact semisimple tensor 2-categories.
\end{Corollary}

Below, we will use the above isomorphisms to exhibit non-trivial classes in $\mathrm{Witt}(\mathbf{Vect}_{\mathds{k}})$. Before doing so, we will establish a slightly weaker criterion for invertibility and examine a pertinent class of examples.

\begin{Lemma}\label{lem:weakerinvertibility}
A compact semisimple tensor 2-category $\mathfrak{C}$ over $\mathds{k}$ is invertible if and only if $\Omega\mathscr{Z}(\mathfrak{C})\simeq \mathbf{Vect}_{\mathds{k}}$.
\end{Lemma}
\begin{proof}
The forward direction is evident.
As for the backward direction, it follows from \cite[Proof of Thm.\ 3.1.4]{D10} and \cite[Lemma 2.2.1]{D9} that taking the Drinfeld center commutes with base change to the algebraic closure. But, if $\mathfrak{D}$ is a fusion 2-category over $\overline{\mathds{k}}$ such that $\Omega\mathscr{Z}(\mathfrak{C})\simeq \mathbf{Vect}_{\overline{\mathds{k}}}$, then we must have $\mathscr{Z}(\mathfrak{C})\simeq \mathbf{2Vect}_{\overline{\mathds{k}}}$ by \cite[Theorem 4.2.2 and Corollary 3.2.3]{D9}.
\end{proof}

\begin{Example}\label{ex:invertibleCSST2C}
Let $\mathbb{K}/\mathds{k}$ be a finite normal extension with Galois group $\Gamma = \mathrm{Gal}(\mathbb{K}/\mathds{k})$. 
For any 4-cocycle representing a class $\pi$ in $H^4(\Gamma;\overline{\mathds{k}}^{\times})$, we can consider the corresponding compact semisimple tensor 2-category $\mathbf{2Vect}_{\mathbb{K}}^{\pi}(\Gamma)$ over $\mathds{k}$ as in example \ref{ex:gradedCSS2C} above.
By construction, note that the connected component of $\mathbf{2Vect}_{\mathbb{K}}^{\pi}(\Gamma)$ with grading $\gamma\in \Gamma$ contains a preferred invertible object $\mathbf{Vect}_{\mathbb{K}}$ that we denote by $\mathbf{V}^{\pi}_{\gamma}$.
We assert that these compact semisimple tensor 2-categories are (Morita) invertible.
Namely, we have $\Omega\mathscr{Z}\big(\mathbf{2Vect}^{\pi}_{\mathbb{K}}(\Gamma)\big)\simeq \mathbf{Vect}_{\mathds{k}}$ as this computation does not involve the pentagonator. The claim therefore follows from Lemma \ref{lem:weakerinvertibility}.

Inspired by \cite[Theorem 3.6]{SS}, we also examine the tensor product of two such compact semisimple tensor 2-categories. More precisely, let $\pi$ and $\varpi$ be two classes in $H^4(\Gamma;\overline{\mathds{k}}^{\times})$. The compact semisimple multitensor 2-category $$\mathfrak{C}:=\mathbf{2Vect}_{\mathbb{K}}^{\pi}(\Gamma)\boxtimes_{\mathds{k}}\mathbf{2Vect}_{\mathbb{K}}^{\varpi}(\Gamma)$$ is invertible, and must therefore be indecomposable. We claim that it is Morita equivalent to $\mathbf{2Vect}_{\mathbb{K}}^{\pi\varpi}(\Gamma)$.
More precisely, we will identify $\mathbf{2Vect}_{\mathbb{K}}^{\pi\varpi}(\Gamma)$ as one of the diagonal blocks in the matrix decomposition of $\mathfrak{C}$ arising from the splitting of its monoidal unit into a direct sum of simple objects.
Our analysis is a pedestrian categorification of \cite[Proposition 3.5]{SS}.
Firstly, we have an isomorphism $$\Omega^2\mathfrak{C} = \mathbb{K}\otimes_{\mathds{k}}\mathbb{K} \cong \prod_{\gamma\in\Gamma}\mathbb{K}$$ of $\mathds{k}$-algebras.
More precisely, the element $a\otimes b$ of $\mathbb{K}\otimes_{\mathds{k}}\mathbb{K}$ is sent to the tuple $(c_{\gamma})_{\gamma\in\Gamma}$ with $c_{\gamma}:=a\gamma(b)$. In particular, the canonical action of $\Gamma\times\Gamma$ on $\mathbb{K}\otimes_{\mathds{k}}\mathbb{K}$ translates to the action of $\Gamma\times\Gamma$ on $\prod_{\gamma\in\Gamma}\mathbb{K}$ given by $(\alpha,\beta)\cdot c_{\gamma} = \alpha(c_{\alpha^{-1}\gamma\beta})$.
We write $I_{\delta}$ for the summand of the monoidal unit $I$ of $\mathfrak{C}$ corresponding to the tuples whose only non-zero entry is $\delta\in\Gamma$.
It follows that, for every triple $\alpha,\beta,\delta\in\Gamma$, the object $$I_{\delta}\Box (\mathbf{V}_{\alpha}^{\pi}\boxtimes \mathbf{V}_{\beta}^{\varpi})$$ in $\mathfrak{C}$ is simple.
The above collection of simple objects contains exactly one simple object in every connected component of $\mathfrak{C}$.
Moreover, we have $$ (\mathbf{V}_{\alpha}^{\pi}\boxtimes \mathbf{V}_{\beta}^{\varpi})\Box I_{\delta}\cong I_{\alpha^{-1}\delta\beta}\Box (\mathbf{V}_{\alpha}^{\pi}\boxtimes \mathbf{V}_{\beta}^{\varpi}).$$ Putting the above observations together, it follows that the monoidal 2-functor $\mathbf{2Vect}_{\mathbb{K}}^{\pi\varpi}(\Gamma)\rightarrow \mathfrak{C}$ induced by $\mathbf{V}_{\delta}^{\pi\varpi}\mapsto I_{\delta^{-1}}\Box (\mathbf{V}_{\delta}^{\pi}\boxtimes \mathbf{V}_{\delta}^{\varpi})$ witnesses the inclusion of a diagonal summand. This concludes the proof of the claim.
\end{Example}

The next result should be compared both with \cite[Theorem 5.9]{SS} and the classical identification by Noether of the Brauer group with the second Galois cohomology of the absolute Galois group.

\begin{Theorem}\label{thm:GaloisH4inclusion}
There is an monomorphism of groups $$H^4(\mathrm{Gal}(\overline{\mathds{k}}/\mathds{k});\overline{\mathds{k}}^{\times})\hookrightarrow \mathrm{Witt}(\mathbf{Vect}_{\overline{\mathds{k}}}).$$
\end{Theorem}
\begin{proof}
By definition, we have $$H^4(\mathrm{Gal}(\overline{\mathds{k}}/\mathds{k});\overline{\mathds{k}}^{\times}) = \mathrm{colim}_{\mathbb{K}}H^4(\mathrm{Gal}(\mathbb{K}/\mathds{k});\mathbb{K}^\times),$$ where the colimit ranges over all finite Galois extensions of $\mathds{k}$, and is taken along the inflations maps $\mathrm{Inf}_{\mathbb{K}}^{\mathbb{L}}:H^4(\mathrm{Gal}(\mathbb{K}/\mathds{k});\mathbb{K}^\times)\rightarrow H^4(\mathrm{Gal}(\mathbb{L}/\mathds{k});\mathbb{L}^\times)$, where $\mathbb{L}/\mathbb{K}$ and $\mathbb{K}/\mathds{k}$ are finite Galois extensions.
In particular, any class in $H^4(\mathrm{Gal}(\overline{\mathds{k}}/\mathds{k});\overline{\mathds{k}}^{\times})$ can be represented by a class $\pi$ in $H^4(\mathrm{Gal}(\mathbb{K}/\mathds{k});\mathbb{K}^\times)$ for some finite Galois extension $\mathbb{K}$ of $\mathds{k}$.
We can consider the corresponding invertible compact semisimple tensor 2-category $\mathbf{2Vect}^{\pi}_{\mathbb{K}}(\mathrm{Gal}(\mathbb{K}/\mathds{k}))$ over $\mathds{k}$.
Additionally, by example \ref{ex:invertibleCSST2C}, the assignment $\pi\mapsto \mathbf{2Vect}^{\pi}_{\mathbb{K}}(\mathrm{Gal}(\mathbb{K}/\mathds{k}))$ defines a group homomorphism $H^4(\mathrm{Gal}(\mathbb{K}/\mathds{k});\mathbb{K}^\times)\rightarrow \mathrm{Witt}(\mathbf{Vect}_{\mathds{k}})$.
We claim that they can be assembled into a group homomorphism $H^4(\mathrm{Gal}(\overline{\mathds{k}}/\mathds{k});\overline{\mathds{k}}^{\times})\rightarrow \mathrm{Witt}(\mathbf{Vect}_{\overline{\mathds{k}}})$.
This follows from the observation that there is an equivalence \begin{equation}\label{eq:inflationbothways}\mathrm{Inf}_{\mathbb{K}}^{\mathbb{L}}\big(\mathbf{2Vect}_{\mathbb{K}}^{\pi}(\mathrm{Gal}(\mathbb{K}/\mathds{k}))\big)\simeq \mathbf{2Vect}_{\mathbb{K}}^{\mathrm{Inf}_{\mathbb{K}}^{\mathbb{L}}(\pi)}(\mathrm{Gal}(\mathbb{L}/\mathds{k}))\end{equation} of compact semisimple tensor 2-categories over $\mathds{k}$, which follows from a straightforward categorification of \cite[Theorem 3.13]{SS}.

It therefore only remains to check that $H^4(\mathrm{Gal}(\overline{\mathds{k}}/\mathds{k});\overline{\mathds{k}}^{\times})\rightarrow \mathrm{Witt}(\mathbf{Vect}_{\mathds{k}})$ is a monomorphism.
Said differently, we have to show that if the class $\pi$ in $H^4(\mathrm{Gal}(\mathbb{K}/\mathds{k});\mathbb{K}^\times)$ is such that $\mathfrak{C}:=\mathbf{2Vect}^{\pi}_{\mathbb{K}}(\mathrm{Gal}(\mathbb{K}/\mathds{k}))$ is Morita equivalent to $\mathbf{2Vect}_{\mathds{k}}$, then there exists a finite Galois extension $\mathbb{L}/\mathbb{K}$ such that $\mathrm{Inf}_{\mathbb{K}}^{\mathbb{L}}(\pi) = \mathrm{triv}$.
Let $\mathfrak{M}$ be an invertible $\mathfrak{C}$-$\mathbf{2Vect}_{\mathds{k}}$-bimodule 2-category.
As the compact semisimple tensor 2-category $\mathbf{2Vect}_{\mathds{k}}$ is connected, we must have $\mathfrak{M}\simeq \mathbf{Mod}(\mathcal{C})$ for some fusion 1-category $\mathcal{C}$ over $\mathds{k}$.
In fact, $\mathfrak{M}$ is also a module 2-category over $\mathbf{2Vect}_{\mathbb{K}}\subseteq \mathfrak{C}$, so that $\mathcal{C}$ is in fact linear over $\mathbb{K}$. Moreover, it follows from Lemma \ref{lem:loopsdualconnected} that $\mathcal{Z}(\mathcal{C},\mathbf{Vect}_{\mathds{k}}) = \mathcal{Z}(\mathcal{C})\simeq \mathbf{Vect}_{\mathbb{K}}$, so that $\mathcal{C}$ is a Morita invertible fusion 1-category over $\mathbb{K}$.
Let $\mathbb{L}/\mathbb{K}$ be a finite Galois extension such that the base change $\mathcal{C}_{\mathbb{L}}$ is Morita equivalent to $\mathbf{Vect}_{\mathbb{L}}$. Such an extension must exist by Lemma \ref{lem:equivalencefinitebasechange} as $\overline{\mathcal{C}}$ is necessarily Morita equivalent to $\mathbf{Vect}_{\overline{\mathds{k}}}$. By construction, the inflation is given by $$\mathrm{Inf}_{\mathbb{K}}^{\mathbb{L}}(\mathfrak{C}) = \mathbf{Bimod}_{\mathfrak{C}}(\mathbf{Vect}_{\mathbb{L}}),$$ where $\mathbf{Vect}_{\mathbb{L}}$ is viewed as a connected rigid algebra in $\mathbf{2Vect}_{\mathbb{K}}\subseteq\mathfrak{C}$.
In particular, we have a canonical invertible bimodule 2-category $\mathbf{Mod}_{\mathfrak{C}}(\mathbf{Vect}_{\mathbb{L}})$ witnessing the Morita equivalence between $\mathbf{2Vect}_{\mathbb{K}}^{\pi}(\mathrm{Gal}(\mathbb{K}/\mathds{k}))$ and its inflation.
It follows that $\mathbf{Mod}_{\mathfrak{M}}(\mathbf{Vect}_{\mathbb{L}})$ is a invertible bimodule 2-category witnessing a Morita equivalence between $\mathrm{Inf}_{\mathbb{K}}^{\mathbb{L}}(\mathfrak{C})$ and $\mathbf{Vect}_{\mathds{k}}$.
But, by unpacking the definitions, we find that $$\mathbf{Mod}_{\mathfrak{M}}(\mathbf{Vect}_{\mathbb{L}})\simeq \mathbf{Mod}(\mathcal{C}_{\mathbb{L}})\simeq \mathbf{2Vect}_{\mathbb{L}}$$ as compact semisimple 2-categories by choice of $\mathbb{L}$.
This implies that the inflation $\mathrm{Inf}_{\mathbb{K}}^{\mathbb{L}}(\mathfrak{C})$ is the Morita dual to $\mathbf{2Vect}_{\mathds{k}}$ with respect to the compact semisimple module 2-category $\mathbf{2Vect}_{\mathbb{L}}$.
Finally, combining equation \eqref{eq:inflationbothways} together with example \ref{ex:MoritatrivialCSS2Cs}, we find that $\mathrm{Inf}_{\mathbb{K}}^{\mathbb{L}}(\pi)=\mathrm{triv}$, thereby concluding the proof.
\end{proof}

In order for the last theorem above to be interesting, we need to exhibit fields for which the fourth cohomology of the absolute Galois group is non-trivial. A large class of examples is provided by totally real number fields, that is, algebraic extensions of $\mathbb{Q}$ for which every embedding into $\mathbb{C}$ factors through $\mathbb{R}$.

\begin{Lemma}\label{lem:GaloisH4}
Let $\mathbb{K}/\mathbb{Q}$ be a totally real algebraic extension. Then, we have $\mathbb{Z}/2\subseteq H^4(\mathrm{Gal}(\overline{\mathbb{Q}}/\mathbb{K}); \overline{\mathbb{Q}}^{\times})$.
\end{Lemma}
\begin{proof}
As $\mathbb{K}$ is totally real, it does not contain any root of unity besides $\pm 1$. In particular, $\mathbb{L}:=\mathbb{K}[\sqrt{-1}]$ is a degree two normal extension of $\mathbb{K}$. We have a semi-direct product decomposition $$\mathrm{Gal}(\overline{\mathbb{Q}}/\mathbb{K})\cong \mathrm{Gal}(\mathbb{L}/\mathbb{K})\ltimes \mathrm{Gal}(\overline{\mathbb{Q}}/\mathbb{L}).$$ As $\mathrm{Gal}(\mathbb{L}/\mathbb{K})\cong\mathbb{Z}/2$, it follows from a straightforward explicit computation that $$H^n(\mathrm{Gal}(\mathbb{L}/\mathbb{K});\overline{\mathbb{Q}}^{\times})\cong \begin{cases}\overline{\mathbb{Q}}^{\times} & \textrm{if}\ n=0,\\0 & \textrm{if}\ n\ \textrm{is odd},\\ \mathbb{Z}/2 & \textrm{if}\ n>0\ \textrm{is even}.\end{cases}$$ This concludes the proof.
\end{proof}

\begin{Remark}\label{rem:exoticexamples}
We suspect that more exotic examples exist. More precisely, let us write $\mathbb{K}:=\mathds{k}(x,y,z)$, a transcendental extension of degree three. It was recorded in \cite[Example 2.8]{SS} that, provided that $\mathds{k}$ has enough roots of unity, then $H^3(\mathrm{Gal}(\overline{\mathbb{K}}/\mathbb{K});\overline{\mathbb{K}}^{\times})$ is non-zero. We imagine that, if $\mathbb{L}:=\mathds{k}(w,x,y,z)$ is a transcendental extension of degree four and $\mathds{k}$ has enough roots of unity, then $H^4(\mathrm{Gal}(\overline{\mathbb{L}}/\mathbb{L});\overline{\mathbb{L}}^{\times})$ is non-zero.
\end{Remark}

\subsection{Witt Spaces and a Conjecture}\label{sub:SpectralSeq}

From the symmetric monoidal 4-category $\mathbf{BrFus}_{\mathds{k}}$ of braided fusion 1-categories over $\mathds{k}$, we can extract the spectrum $\mathcal{W}itt(\mathbf{Vect}_{\mathds{k}}):=\mathbf{BrFus}_{\mathds{k}}^{\times}$ consisting of invertible objects and invertible morphisms. The homotopy groups of $\mathcal{W}itt(\mathbf{Vect}_{\mathds{k}})$ are given by: $$\begin{tabular}{|c|c|c|c|c|}
\hline
$\pi_0$ & $\pi_1$ & $\pi_2$ & $\pi_3$ & $\pi_4$\\
\hline \\[-1em]
$\mathrm{Witt}(\mathbf{Vect}_{\mathds{k}})$ & $\mathrm{MoInv}(\mathbf{Vect}_{\mathds{k}})$ & $\mathrm{Br}(\mathds{k})$ & $0$ & $\mathds{k}^{\times}$\\
\hline
\end{tabular}$$
Note that those groups are precisely those recorded in table \ref{tab:ourbelovedgroups}. More generally, using the techniques of \cite{BJS}, one can associate to any symmetric (multi)fusion 1-category $\mathcal{E}$ over $\mathds{k}$ a symmetric monoidal 4-category $\mathbf{BrFus}(\mathcal{E})$ of braided fusion 1-categories enriched over $\mathcal{E}$. We write $\mathcal{W}itt(\mathcal{E})$ for the corresponding spectrum of invertible objects and invertible morphisms.

Let us write $\mathbf{SMF1C}^{\mathrm{dom,faith}}_{\overline{\mathds{k}}}$ for the $(2,1)$-category of symmetric (multi)fusion 1-categories over the algebraically closed field $\overline{\mathds{k}}$ with dominant faithful symmetric tensor functors.
It was shown in \cite[Proposition 1.13]{DHJFNPPRY} that there is a functor $\mathcal{W}itt:\mathbf{SMF1C}^{\mathrm{dom,faith}}_{\overline{\mathds{k}}}\rightarrow \mathbf{Spaces}$ sending a symmetric (multi)fusion 1-category $\overline{\mathcal{E}}$ over $\overline{\mathds{k}}$ to the Witt space $\mathcal{W}itt(\overline{\mathcal{E}})$. Moreover, it is expected that the fact that every fusion 2-categories is Morita equivalent to a connected one can be used to show that the functor $\mathcal{W}itt$ commutes with (homotopy) limits of finite groups.

Thanks to Theorem \ref{maintheorem}, we suspect that the above conjecture can be generalized to an arbitrary field $\mathds{k}$.
More precisely, let us write $\widehat{\mathbf{SMF1C}}^{\mathrm{dom,faith}}_{\mathds{k}}$ for the $(2,1)$-category of extended symmetric (multi)fusion 1-category over $\mathds{k}$ and faithful dominant symmetric tensor functors.
By an extended symmetric fusion 1-category, we mean a symmetric monoidal 1-category $\mathcal{E}$ such that $\Omega\mathcal{E}$ is an algebraic (but not necessarily finite) extension of $\mathds{k}$ and such that $\mathcal{E}$ is symmetric fusion over $\Omega\mathcal{E}$.
An extended symmetric multifusion 1-category is a finite direct sum of extended symmetric fusion 1-categories.
Thus, for any symmetric (multi)fusion 1-category $\mathcal{E}$ over $\mathds{k}$, its base change $\overline{\mathcal{E}}$ to the algebraic closure is also an object in $\widehat{\mathbf{SMF1C}}^{\mathrm{dom,faith}}_{\mathds{k}}$.

\begin{Conjecture}
There is a functor $\mathcal{W}itt:\widehat{\mathbf{SMF1C}}^{\mathrm{dom,faith}}_{\mathds{k}}\rightarrow \mathbf{Spaces}$ from the $(2,1)$-category of extended symmetric (multi)fusion 1-categories over $\mathds{k}$ to the $(\infty,1)$-category of spaces.
This functor preserves homotopy limits of profinite groups.
\end{Conjecture}

The above conjecture would in particular make rigorous the discussion of \cite[Section 6]{SS}, but also yields a computational method to study the groups $\mathrm{Witt}(\mathcal{E}) = \pi_0(\mathcal{W}itt(\mathcal{E}))$ via Lemma \ref{lem:finitedescent}. As a case study, we now discuss the spectral sequence associated to the action of the absolute Galois group $\Gamma = \mathrm{Gal}(\overline{\mathds{k}}/\mathds{k})$ on $\overline{\mathds{k}}$. The above conjecture gives a homotopy fixed point spectral sequence $$E_2^{p,q} = H^p(\Gamma, \pi_q\mathcal{W}itt(\mathbf{Vect}_{\overline{\mathds{k}}}))\Rightarrow \pi_{q-p}\mathcal{W}itt(\mathbf{Vect}_{\mathds{k}}),$$ whose $E_2$-page is given by:
$$\renewcommand{\arraystretch}{1.4}
    \begin{array}{c|ccccccc}
      q \\
      \\
     4 & \mathds{k}^{\times} & 0 & H^2(\Gamma;\overline{\mathds{k}}^{\times}) & H^3(\Gamma;\overline{\mathds{k}}^{\times}) & H^4(\Gamma;\overline{\mathds{k}}^{\times}) \\
     & \\
     0 & \mathrm{Witt}(\mathbf{Vect}_{\overline{\mathds{k}}})^{\Gamma}\\
     \hline 
     & 0 & 1 & 2 & 3 & 4 & 5 & \quad p\,.
\end{array}$$
\noindent In particular, the only potentially interesting differential is an edge differential $$d_4:E^{0,0}_2\cong \mathrm{Witt}(\mathbf{Vect}_{\overline{\mathds{k}}})^{\Gamma}\rightarrow E^{5,4}_4\cong H^5(\Gamma;\overline{\mathds{k}}^{\times}).$$ Consequently, as already observed in \cite{SS}, we obtain isomorphisms $$\mathrm{Br}(\mathds{k})\cong H^2(\Gamma;\overline{\mathds{k}}^{\times}),\, \mathrm{MoInv}(\mathbf{Vect}_{\mathds{k}})\cong H^3(\Gamma;\overline{\mathds{k}}^{\times}),$$ but also an exact sequence \begin{equation}0\rightarrow H^4(\Gamma;\overline{\mathds{k}}^{\times})\rightarrow \mathrm{Witt}(\mathbf{Vect}_{\mathds{k}})\rightarrow \mathrm{Witt}(\mathbf{Vect}_{\overline{\mathds{k}}})^{\Gamma}\rightarrow H^5(\Gamma;\overline{\mathds{k}}^{\times}),\end{equation}
which provides a different perspective on Theorem \ref{thm:GaloisH4inclusion}.

\begin{Example}
In the case $\mathds{k} = \mathbb{R}$, we therefore have a short exact sequence $$0\rightarrow \mathbb{Z}/2\rightarrow \mathrm{Witt}(\mathbf{Vect}_{\mathbb{R}})\rightarrow \mathrm{Witt}(\mathbf{Vect}_{\mathbb{C}})^{\Gamma}\rightarrow 0,$$ as $H^5(\Gamma;\mathbb{C}^{\times})=0$ by the proof of Lemma \ref{lem:GaloisH4}. On the other hand, the structure of the abelian group $\mathrm{Witt}(\mathbf{Vect}_{\mathbb{C}})^{\Gamma}$ is difficult to identify because we do not know explicit generators for the abelian group $\mathrm{Witt}(\mathbf{Vect}_{\mathbb{C}})$.
\end{Example}

\begin{Remark}
The above conjecture, or more precisely the spectral sequences it provides, also offers a different perspective on Theorem \ref{thm:Picardgrouptorsion}. Succinctly, let $\mathcal{E}$ be a symmetric fusion 1-category over $\mathds{k}$. Then, $\overline{\mathcal{E}}$ is a symmetric fusion 1-category over $\overline{\mathds{k}}$ and the homotopy groups $\pi_q\mathcal{W}itt(\overline{\mathcal{E}})$ are torsion for $q=1,2,3$. This follows, for instance, by combining the above conjecture with Deligne's Theorem \cite{De}. But then, as $\mathcal{E}\simeq \overline{\mathcal{E}}^{\Gamma}$, the corresponding spectral sequence shows that $\mathrm{Pic}(\mathcal{E}) = \pi_2\mathcal{W}itt(\mathcal{E})$ is a torsion abelian group. In fact, the proof of Theorem \ref{thm:Picardgrouptorsion} is essentially a purely algebraic version of this homotopy theoretic argument.
In addition, we also find that $\mathrm{MoInv}(\mathcal{E})=\pi_1\mathcal{W}itt(\mathcal{E})$ is torsion. On the other hand, we can not expect $\mathrm{Witt}(\mathcal{E})$ to be torsion because $\mathrm{Witt}(\overline{\mathcal{E}})$ is not.
\end{Remark}

\bibliography{bibliography.bib}

\end{document}